\documentclass[a4paper,11pt,twoside]{article}
\usepackage{amsmath,amssymb}
\usepackage{float}
\usepackage[
sort&compress]{natbib}

\usepackage[text={6.5in,9.1in},centering]{geometry}

\usepackage{theorem}

\newtheorem{theorem}{Theorem}[section]
\newtheorem{proposition}{Proposition}[section]
\newtheorem{corollary}{Corollary}[section]
\newtheorem{lemma}{Lemma}[section]

{\theorembodyfont{\upshape}
\newtheorem{remark}{Remark}[section]
\newtheorem{example}{Example}[section]
}

\newcommand{\Proof}{\textbf{Proof. }}
\newcommand{\Proofof}[1]{\textbf{Proof of #1. }}

\newcommand{\CC}{\mathbb{C}}
\newcommand{\RR}{\mathbb{R}}
\newcommand{\NN}{\mathbb{N}}
\newcommand{\ZZ}{\mathbb{Z}}
\newcommand{\Zpl}{\ZZ_+}

\newcommand{\A}{\mathcal{A}}
\newcommand{\F}{\mathcal{F}}
\newcommand{\M}{\mathcal{M}}

\newcommand{\X}{\mathfrak{X}}

\newcommand{\cc}{{\mathrm{c}}}
\newcommand{\dd}{{\mathrm{d}}}
\newcommand{\ee}{{\mathrm e}}
\newcommand{\ii}{{\mathrm i}}
\newcommand{\pii}{{\mathrm{\pi}}}

\newcommand{\ab}[1]{|#1|}                   

\newcommand{\Ab}[1]{\Big|#1\Big|}           

\newcommand{\bessel}{\mathrm{Bessel}}       

\newcommand{\bydef}{=:}                     

\newcommand{\card}[1]{|#1|}                 

\newcommand{\ceil}[1]{{\lceil #1\rceil}}

\newcommand{\coeff}[2]{\mathrm{Coeff}(#1,#2)}

\newcommand{\Coeff}[2]{\mathrm{Coeff}\Big(#1,#2\Big)}

\newcommand{\defby}{:=}                     

\newcommand{\dirac}{I}                      

\newcommand{\Expect}{\mathrm{E}}            

\newcommand{\expo}[1]{\exp(#1)}             
\newcommand{\Expo}[1]{\exp\Big(#1\Big)}     

\newcommand{\floor}[1]{{\lfloor #1\rfloor}}
\newcommand{\Floor}[1]{{\Big\lfloor #1\Big\rfloor}}

\renewcommand{\geq}{\geqslant}

\newcommand{\indicator}{\mathbf{1}}         

\newcommand{\Kr}{\mathrm{Kraw}}             

\renewcommand{\leq}{\leqslant}

\newcommand{\muc}{\mathrm{Mult}}      

\newcommand{\norm}[1]{\big\|#1\big\|}       

\newcommand{\Norm}[1]{\Big\|#1\Big\|}       

\newcommand{\Prob}{\mathrm{P}}              

\newcommand{\qed}{$\square$}                

\newcommand{\restr}[2]{#1|_{#2}^{}}   
 
\newcommand{\set}[1]{\underline{#1}}

\newcommand{\setn}[1]{\underline{#1}_0} 
 
\newcommand{\tvm}[1]{|#1|}                  

\newcommand{\unitvec}[1]{e_{#1}^{}}

\newcommand{\vecsum}[1]{|#1|}

\newcommand{\binomial}[2]{\genfrac{(}{)}{0pt}{}{#1}{#2}}

\newcommand{\newatop}[2]{\genfrac{}{}{0pt}{}{#1}{#2}}

\scrollmode
\begin{document}
\title{
Closeness of convolutions of probability measures}
\author{Bero Roos\medskip\\
Department of Mathematics,
University of Leicester,
University Road, \\
Leicester LE1 7RH,
United Kingdom. 
E-mail: b.roos@leicester.ac.uk\bigskip\\
Running title: Closeness of convolutions
}
\date{Revised version
}
\maketitle
\begin{abstract}
We derive new explicit bounds for
the total variation distance between two convolution products of
$n\in\NN$ probability distributions, one
of which having identical convolution factors. Approximations by
finite signed measures of arbitrary order are considered as well. We
are interested in bounds with magic factors, i.e.\ roughly
speaking~$n$ also appears  in the denominator. Special emphasis is
given to the approximation by the $n$-fold convolution of the
arithmetic mean of the distributions under consideration. As an
application, we consider the multinomial approximation  of the
generalized multinomial distribution. It turns out that here  the 
order of some bounds given in \citet{roo01} and \citet{loh92} can 
significantly be improved. In particular, it follows that a dimension
factor can be dropped. Moreover, better accuracy is achieved in the
context of symmetric distributions with finite support. 
In the course of proof, we use a basic Banach algebra technique
for measures on a measurable Abelian group.
Though this method was already used by \citet{lec60}, 
our central arguments seem to be new.
We also derive new smoothness bounds for convolutions of probability 
distributions, which might be of independent interest.\medskip \\
{\emph{Keywords: }\small    
Convolutions, explicit constants, generalized multinomial 
distribution, multivariate Krawtchouk polynomials, 
magic factor, multinomial approximation, 
signed measures, total variation distance.
\medskip\\
}
{\small MSC 2000 Subject Classification:
Primary 60F05;   
secondary 60G50, 
62E17.           
}
\end{abstract}
\section{Introduction}
\subsection{Aim of the paper}
Approximations of distributions of sums of independent random 
variables are needed in nearly all branches of probability theory and
statistics. Many results for normal and compound Poisson
approximations are nowadays available.  However, if the distributions
of the summands are similar to each other,  much better accuracy can
be achieved using identical convolutions of a certain  distribution.
In the present paper, we give total variation bounds for the accuracy
of such approximations in a general framework,  i.e.~for 
probability distributions on a measurable Abelian group.  We also
consider higher order approximations by finite signed  measures.  All
bounds contain magic factors, i.e.~roughly speaking~$n$ appears in the
denominator. As a consequence, this enables us to derive
multidimensional results, some of which improve the order of bounds
obtained in \citet{roo01} and \citet{loh92}. It should be mentioned 
that Loh used Stein's method  in a more general situation of dependent
random variables. However, it seems to be unclear, whether Stein's
method can be  used to reproduce the results of the present paper.
Furthermore, it turns out that our bounds have a better order
in the case of symmetric probability distributions with finite 
support. Our proofs are based on a combination of some Banach algebra
related techniques, which in principle were used by \citet{lec60}. 
On the other hand, the core arguments given in
Sections~\ref{aux} and~\ref{genla} seem to be new.
Further, the smoothness estimates for convolutions of  
probability distributions in Section~\ref{aux} might be of  
independent interest; for instance, see (\ref{eq18798}) and 
(\ref{eq1843}).

We note that, at the beginning of our investigation, we tried to 
improve one of the central results of \citet{roo01}, see (\ref{eq1}) 
and discussion thereafter. But unfortunately we were 
not able to use the multidimensional expansion of that paper for any 
substantial improvement. Surprisingly it turned out that it is
better to forget the dimension, so to speak, and to use the 
properties of measures on a measurable Abelian group. 
This should explain, why we use this somewhat abstract approach.

The paper is structured as follows: The following two subsections 
are devoted to the notation and a review of known results. In
Section~\ref{mainres}, we present and discuss our main
results. To get a first impression of the results of this paper,
the reader may consult (\ref{eq188}), (\ref{eq188b}), 
and (\ref{eq7884}). In Section~\ref{numerics}, we give some numerical
examples. The proofs are contained in Section~\ref{proofs}.
\subsection{Notation}
Let $(\X,+,\A)$ be a measurable Abelian group,  that is, $(\X,+)$ is a
commutative group with identity element~$0$ and $\A$ is a
$\sigma$-algebra of subsets of $\X$  such that the mapping 
$(x,y)\mapsto x-y$ from $(\X\times\X,\A\otimes\A)$ to $(\X,\A)$ is
measurable. We note that it is more convenient to formulate our
results in terms of distributions or signed measures rather than in
terms of random variables. Let $\F$ (resp.~$\M$) be the set of all 
probability distributions (resp.~finite signed measures) on 
$(\X,\A)$. Products and powers of finite signed measures in $\M$ are
defined in the convolution  sense, that is, for $V,W\in\M$ and
$A\in\A$, we write $VW(A)=\int_\X V(A-x)\,\dd W(x)$. Empty products
and powers of  signed measures in $\M$ are understood to be
$\dirac\defby\dirac_0$, where  $\dirac_x$  is the Dirac measure at 
point $x\in\X$. Let $V=V^{+}-V^{-}$
denote the Hahn-Jordan decomposition of~$V\in\M$ and let
$\tvm{V}=V^{+}+V^{-}$ be its total variation measure. The total
variation norm of $V$ is defined by  $\norm{V}=\tvm{V}(\X)$. 
We note that, in the literature, often the total variation distance 
$\sup_{A\in\A}\ab{F(A)-G(A)}=\frac{1}{2}\norm{F-G}$ 
between $F,G\in\F$ is used. In this paper, however, all distances
will be given only in the total variation norm. With the 
usual operations of real scalar multiplication, addition,  together
with convolution and the total variation norm, $\M$ is a real
commutative Banach algebra with unity $\dirac$.
For $V\in\M$ and a power series
$g(z)=\sum_{m=0}^\infty a_mz^m$, $(a_m\in\RR)$ 
converging absolutely for each complex $z\in\CC$ with 
$\ab{z}\leq\norm{V}$, we define $g(V)=\sum_{m=0}^\infty a_mV^m$. 
The above assumptions imply that the limit exists and is an element 
of the Banach algebra $\M$. On the other hand, the definition of 
$g(V)$ can also be understood setwise.
The exponential of $V\in\M$ is defined by the finite signed measure
\[\ee^V=\expo{V}\,=\,\sum_{m=0}^\infty\frac{1}{m!}V^m\in\M.\] 
We note that $\ee^V$ is not necessarily a non-negative measure.
Further, $\expo{t(F-\dirac)}$ is the compound Poisson
distribution with parameters $t\in[0,\infty)$ and $F\in\F$.  
If $F$ and $G$ are non-negative measures on $(\X,\A)$ and $F$ is 
absolutely continuous with respect to $G$, we write $F\ll G$.
For  $F\in\F$ and $A\in\A$, $\restr{F}{A}$ is the restriction of $F$
to the set $A$. The complement of $A\in\A$ is denoted by $A^\cc$.
Set $\set{0}=\emptyset$ and $\set{n}=\{1,\dots,n\}$
for $n\in\NN=\{1,2,\dots\}$; further, for  $n\in\Zpl=\NN\cup\{0\}$, 
set $\setn{n}=\{0,\dots,n\}$. For a set $J$, let $\card{J}$ be the 
number of its elements. For $x\in\RR$, let 
$\floor{x}=\sup\{n\in\ZZ\,|\,n\leq x\}$ and
$\ceil{x}=\inf\{n\in\ZZ\,|\,n\geq x\}$. 
Always, let $0^0=1$, $1/0=\infty$,  and, for $k\in\ZZ$, 
$\sum_{m=k}^{k-1}=0$ be the empty sum and
$\prod_{m=k}^{k-1}=1$ the empty product.
For $a\in\CC$ and $b\in\Zpl$, let 
$\binomial{a}{b}=\prod_{m=1}^b(a-m+1)/m$.
For $a,b\in\RR$, set $a\wedge b=\min\{a,b\}$. 
\subsection{Known results}\label{known}
We first discuss some important results for discrete distributions 
on $\X=\RR^d$, $(d\in\NN)$ with the usual addition. Let 
\begin{equation}\label{eq1199} 
H_0=\dirac,\quad H_r=\dirac_{\unitvec{r}},\quad(r\in\set{d}),\quad 
F_j=\sum_{r=0}^d p_{j,r}H_r,\quad (j\in\set{n},\;n\in\NN),\quad
\overline{F}=\sum_{r=0}^d\overline{p}_rH_r,
\end{equation}
where, for $r\in\setn{d}$, 
$p_{j,r}\in[0,1]$ with $\sum_{r=0}^dp_{j,r}=1$, 
$\overline{p}_r=n^{-1}\sum_{j=1}^np_{j,r}>0$, 
and $\unitvec{r}\in\RR^d$, ($r\neq0$) is the vector with $1$ at 
position $r$ 
and $0$ otherwise.

In the case $d=1$, \citet[Theorem~1 and Lemma~2]{ehm91} proved with 
the help of Stein's method that the total variation  distance between 
the Bernoulli convolution $\prod_{j=1}^nF_j$ and the binomial 
law~$\overline{F}^n$ can be estimated by
\begin{equation}\label{Ehm}
\frac{\gamma_2}{62}\min\Big\{1,\frac{1}{n\overline{p}_1
\overline{p}_0}\Big\}\leq\Norm{\prod_{j=1}^n F_j-\overline{F}^n}
\leq2\gamma_2\min\Big\{1,\frac{1}{n\overline{p}_1\overline{p}_0}
\Big\},
\end{equation}
where $\gamma_k=\sum_{j=1}^n(\overline{p}_1-p_{j,1})^k$, $(k\in\NN)$.
Here, the
estimates depend on the behavior of the so-called magic  factor
$(n\overline{p}_1\overline{p}_0)^{-1}$ 
(cf.\ Introduction in \citet{BHJ92}), and on the closeness of 
all~$p_{j,1}$, $(j\in\set{n})$, which is reflected by~$\gamma_2$.  
In Theorem~3 of \citet{roo00}, a Krawtchouk expansion was used to 
show that an absolute constant $C>0$ exists such that, if 
$\gamma_2>0$, then
\[\Ab{\Norm{\prod_{j=1}^nF_j-\overline{F}^n}
-\theta\sqrt{\frac{2}{\pii\ee}}}\leq C\,\theta\,
\min\Big\{1,\frac{\ab{\gamma_3}}{\gamma_2\sqrt{n\overline{p}_1
\overline{p}_0}}+\frac{1}{n\overline{p}_1\overline{p}_0}
+\theta\Big\},\qquad 
\Big(\theta=\frac{\gamma_2}{n\overline{p}_1\overline{p}_0}\Big).\]
For example, it easily follows that 
$\norm{\prod_{j=1}^nF_j-\overline{F}^n}
\sim\sqrt{2/(\pii\ee)}\,\theta$ as
$\theta\to0$ and $n\overline{p}_1\overline{p}_0\to\infty$.
 Here,~$\sim$ means that the quotient of both sides tends to one.
Further results in this and a more general context can be found in 
\citet{CR06} and the papers cited there.

The multivariate case $d\in\NN$ was investigated by \citet{loh92}
using Stein's method. He gave an estimate for the closeness between 
the generalized  multinomial distribution  $\prod_{j=1}^nF_j$ and the 
multinomial distribution $\overline{F}^n$. This bound contains 
certain functions $C_1,C_2>0$ of $\overline{p}_r$, $(r\in\setn{d})$,
which can be estimated from above by absolute constants, if  all
$\overline{p}_r$'s are uniformly bounded away from~$0$ and~$1$.  In
his Theorem~5, he showed that, if $n\geq2$ and 
$\max\{C_1n^{-1/2},\,C_2[2(n-1)]^{-1}\}\leq1$, then 
\begin{equation}\label{loh} 
\Norm{\prod_{j=1}^nF_j-\overline{F}^n}\leq
2\sum_{j=1}^n\sum_{0\leq r_1<r_2\leq d}
\ab{p_{j,r_1}\overline{p}_{r_2}-p_{j,r_2}\overline{p}_{r_1}}\,
\varepsilon_{r_1,r_2},
\end{equation}
where
\[\varepsilon_{r_1,r_2}=
\frac{C_2}{n-1}\ln\Big(\frac{2(n-1)}{C_2}\Big)
+\Big(\frac{C_2}{2(n-1)}\Big)^2+
\frac{2C_1}{\sqrt{n}}\min\Big\{\prod_{i=1}^n
(1-p_{i,r_1}),\,\prod_{i=1}^n(1-p_{i,r_2})\Big\}.\]
The quantities $C_1,C_2$ can be given explicitly as 
\begin{equation}\label{eq6341}
C_1=\sup_{0\leq r_1< r_2\leq d}[\widetilde{C}_1(r_1,r_2)
\wedge\widetilde{C}_1(r_2,r_1)],\quad 
C_2=\sup_{\newatop{\scriptstyle 0\leq r_1,r_2,r_3\leq d:}{
\scriptstyle r_2\neq r_1,\,r_3\neq r_1}}
\widetilde{C}_2(r_1,r_2,r_3),
\end{equation}
where, for $r_1,r_2,r_3\in\setn{d}$, 
\begin{eqnarray*}
\widetilde{C}_1(r_1,r_2)&=&\Big(\frac{2}{\overline{p}_{r_1}}
+\frac{3}{\overline{p}_{r_2}}
+\frac{1}{\ee\overline{p}_{r_2}(1-\overline{p}_{r_2})}\Big)^{1/2}
+\Big(\frac{1}{2\ee\overline{p}_{r_2}(1-\overline{p}_{r_2})^2}\Big)^{1/2},
\\
\widetilde{C}_2(r_1,r_2,r_3)&=&
\left\{\begin{array}{ll}
\frac{2}{\overline{p}_{r_1}}+\frac{2}{\overline{p}_{r_2}},&\mbox{if }
r_2=r_3,\\
\Big[\frac{1}{\overline{p}_{r_1}}+
\frac{2}{\ee\overline{p}_{r_2}(1-\overline{p}_{r_2})}
+\frac{2}{\ee\overline{p}_{r_3}(1-\overline{p}_{r_3})}&\\
{}+\Big(\frac{2}{\overline{p}_{r_1}}+\frac{3}{\overline{p}_{r_2}}
+\frac{1}{\ee\overline{p}_{r_2}(1-\overline{p}_{r_2})}
\Big)^{1/2}\Big(\frac{2}{\overline{p}_{r_1}}+\frac{3}{\overline{p}_{r_3}}
+\frac{1}{\ee\overline{p}_{r_3}(1-\overline{p}_{r_3})}
\Big)^{1/2}\Big],&\mbox{if } r_2\neq r_3.
\end{array}\right.
\end{eqnarray*}
If $d=1$, then it follows from Ehm's result and the equality
$\sum_{0\leq r_1<r_2\leq d}\ab{p_{j,r_1}
\overline{p}_{r_2}-p_{j,r_2}\overline{p}_{r_1}}
=\ab{\overline{p}_1-p_{j,1}}$, 
that Loh's bound is not of the best  possible order, because of the 
exponent of $\ab{\overline{p}_1-p_{j,1}}$ and the logarithmic term. 
It turned out that a bound better  than~(\ref{loh}) can be given 
using a multivariate Krawtchouk expansion,
see \citet[Theorem 2, Corollary~1]{roo01}. Indeed,
\begin{equation}\label{eq1} 
\Norm{\prod_{j=1}^nF_j-\overline{F}^n}
\leq C_3\Big(\sum_{r=1}^d\sqrt{\delta(r)}\Big)^2,
\end{equation}
where $C_3=\frac{\ee}{2-\sqrt{3}}\leq10.15$ and 
\[\delta(r)=\sum_{j=1}^n(\overline{p}_r-p_{j,r})^2\min\Big\{
\frac{4}{\ee},\,\frac{1}{n\overline{p}_r\overline{p}_0}\Big\},
\quad (r\in\set{d}).\]
A sometimes more precise bound is 
\begin{equation}\label{eq1b} 
\Norm{\prod_{j=1}^nF_j-\overline{F}^n}
\leq \frac{(\sum_{r=1}^d\sqrt{\ee\,\delta(r)})^2}{
1-\sum_{r=1}^d\sqrt{\ee\,\delta(r)}},\quad\mbox{if} \quad
\sum_{r=1}^d\sqrt{\ee\,\delta(r)}<1.
\end{equation}
In contrast to~(\ref{loh}),  for
$d=1$,  the bounds in~(\ref{eq1}) and~(\ref{Ehm}) have the same
order.  In the general case, from~(\ref{eq1}) and Cauchy's
inequality,  it follows that
$\norm{\prod_{j=1}^nF_j-\overline{F}^n}\leq 
C_3d\sum_{r=1}^d\delta(r)$. 
We note that this estimate is of the same accuracy
as~(\ref{eq1}) when the $\delta(r)$, $(r\in\set{d})$ are of similar
magnitude. In view of this bound, one might wonder, whether the 
dimension factor $d$ can be dropped.
However, as shown in \citet[Remark 2 after Proposition~2]{roo01},
this is not generally possible. But if we concentrate on the 
estimate with the magic factors, i.e.\ 
\begin{equation}\label{eq3590} 
\Norm{\prod_{j=1}^nF_j-\overline{F}^n}\leq 
C_3\,d\,\sum_{r=1}^d
\sum_{j=1}^n\frac{(\overline{p}_r-p_{j,r})^2}{n\overline{p}_r
\overline{p}_0},
\end{equation}
the more general results of this paper imply that $C_3\,d$ can indeed
be replaced by the constant $21.88$, see Example~\ref{ex115} below. 
It should be mentioned that  here the $H_r$, $(r\in\setn{d})$
need not just be the Dirac measures as in~(\ref{eq1199}).
\section{Main results}\label{mainres}
In what follows, we present bounds which are small when 
the $F_j\in\F$, $(j\in\set{n})$ are close or when $n$ is large and 
the $F_j$ are not too different. Our first result is the following.
\begin{theorem}\label{thm01}
Let $n\in\NN$, $F_1,\dots,F_n,G\in\F$, 
$F_0=\overline{F}=\frac{1}{n}\sum_{j=1}^nF_j$,
\[V_k=\sum_{J\subseteq\,\set{n}:\;\card{J}=k}
\prod_{j\in J}(F_j-G),\quad(k\in\setn{n}),\qquad 
W_\ell=\sum_{k=0}^\ell V_k\,G^{n-k}, \quad (\ell\in\setn{n}).\]
For $j,k\in\set{n}$ and $m\in\NN$, set
$M_{j,k}=(F_j-G)G^{\floor{(n-k)/k}}$,
$\nu_{k,m}^{}=\sum_{j=1}^n\norm{M_{j,k}}^m$,
and $\widetilde{\nu}_k^{}=\norm{\sum_{j=1}^nM_{j,k}}$. Set 
\begin{eqnarray*} 
&\displaystyle 
\eta_{\ell,\alpha}^{}=\max_{k\,\in\,\set{n}\setminus\set{\ell}}
\Big[\frac{1}{k^{1+\alpha}}\Big(\frac{\widetilde{\nu}_k^2}{4c_1^{}}+
\nu_{k,2}^{}\Big)\Big],\quad (\ell\in\setn{n},\,\alpha\in[0,\infty)),
\qquad \eta_\ell^{}=\eta_{\ell,0}^{},&\\
&\displaystyle 
c_1^{}=\sup_{x\in(0,\infty)}\Big[\frac{\ln(2-(1-x)\ee^x)}{x^2}\Big]
=0.694025\ldots.&
\end{eqnarray*}
\begin{itemize}
\item[\textup{(a)}] 
Let $\alpha\in[0,\infty)$, $\ell\in\setn{n}$, and 
$\beta=\ceil{\alpha(\ell+1)/2}$. 
If $\eta_{\ell}^{}<(2\ee\,c_1^{})^{-1}$, then
\begin{equation}\label{eq384}
\Norm{\prod_{j=1}^nF_j-W_\ell}\leq
(\ell+1)^{\beta}\,\beta!\,
\frac{(2\ee\,c_1^{}\,\eta_{\ell,\alpha}^{})^{(\ell+1)/2}}{
(1-\sqrt{2\ee\,c_1^{}\,\eta_{\ell}^{}})^{\beta+1}}.
\end{equation}
In particular, for $\alpha=0$, we have
\begin{equation}\label{eq384b}
\Norm{\prod_{j=1}^nF_j-W_\ell}\leq
\frac{(2\ee\,c_1^{}\,\eta_{\ell}^{})^{(\ell+1)/2}}{
1-\sqrt{2\ee\,c_1^{}\,\eta_{\ell}^{}}}.
\end{equation}
\item[\textup{(b)}] 
Assume that, for each $j\in\setn{n}$, 
$B_j\in\A$ exists such that $\restr{F_j}{B_j^\cc}\ll G$ and let 
$f_j$ denote a Radon-Nikodym density of $\restr{F_j}{B_j^\cc}$
with respect to $G$. For $\ell\in\setn{n}$, we then have 
{\arraycolsep3pt
\begin{eqnarray}
\eta_{\ell}^{}
&\leq&\frac{1}{4c_1^{}}\Big[\frac{n}{\sqrt{\ell+1}}
\Big(\norm{\restr{(\overline{F}-G)}{B_0}}+
\ab{(\overline{F}-G)(B_0)}\Big)
+\sqrt{2n}\Big(\int_{B^\cc_0}
(f_0-1)^2\,\dd G\Big)^{1/2}\Big]^2\nonumber\\
&&{}+\sum_{j=1}^n
\Big[\frac{1}{\sqrt{\ell+1}}\Big(\norm{\restr{(F_j-G)}{B_j}}
+\ab{(F_j-G)(B_j)}\Big)
+\sqrt{\frac{2}{n}}\Big(\int_{B^\cc_j}(f_j-1)^2\,\dd G
\Big)^{1/2}\Big]^2.\hspace{5mm}\label{eq512}
\end{eqnarray}}%
\end{itemize} 
\end{theorem} 

We note that, if $G=\overline{F}$, then $\widetilde{\nu}_k=0$ and 
$\eta_{\ell,\alpha}^{}$ simplifies to 
$\eta_{\ell,\alpha}^{}=\max_{k\in\set{n}\setminus\set{\ell}}
\frac{\nu_{k,2}^{}}{k^{1+\alpha}}$. 
One might ask why we gave the complicated estimate (\ref{eq384}).
However, it turns out that, in special situations, the order of 
$\eta_{\ell,\alpha}^{}$ for $\alpha>0$ can be much better than that of
$\eta_{\ell}^{}$. See Proposition~\ref{prop256} below
involving a bound for $\eta_{\ell,1}$ instead of just the 
estimate~(\ref{eq512}).
Further, the reason why we formulated Theorem~\ref{thm01}
in its present general form without the assumption that 
$G=\overline{F}$ is given with Lemma~\ref{la1a} and 
Example~\ref{ex021} below.

Let us first discuss the simple case when $\alpha=0$.
\begin{remark}
Let the assumptions of Theorem~\ref{thm01} hold.
In what follows, whenever we consider $V_k$ or $W_k$ for a
specified number $k\in\Zpl$, we assume that $k\leq n$.
\begin{itemize}
\item[\textup{(a)}] 
For $k\in\NN$, let $\Gamma_k=\sum_{j=1}^n(G-F_j)^k$.
We have $V_0=\dirac$, $V_1=n(\overline{F}-G)$, 
$V_2=\frac{1}{2}(n^2(\overline{F}-G)^2-\Gamma_2)$, 
and, similarly as in \citet[formula~(10)]{roo00}, it can be shown that
\[V_k=-\frac{1}{k}\sum_{j=0}^{k-1}V_j\Gamma_{k-j},
\quad (k\in\set{n}).\]
This formula can easily be used to evaluate the signed measures 
$W_\ell$ for a given $\ell$. In particular, we have $W_0=G^n$ and
\[W_1=G^n+n(\overline{F}-G)G^{n-1},\quad 
W_2=G^n+n(\overline{F}-G)G^{n-1}
+\frac{1}{2}(n^2(\overline{F}-G)^2-\Gamma_2)G^{n-2}.\]
\item[\textup{(b)}] 
In the important case $G=\overline{F}$, the formulas above 
become somewhat simpler. Here, we derive
\begin{eqnarray}
V_1&=&0,\quad 
V_2=-\frac{1}{2}\Gamma_2,\quad
V_3=-\frac{1}{3}\Gamma_3,\quad
V_4=\frac{1}{8}\Gamma_2^2-\frac{1}{4}\Gamma_4, \quad
\label{eq612}\\
V_5&=&\frac{1}{6}\Gamma_2\Gamma_3-\frac{1}{5}\Gamma_5, \quad 
V_6=-\frac{1}{48}\,\Gamma_2^3
+\frac{1}{8}\,\Gamma_2\Gamma_4+\frac{1}{18}\,\Gamma_3^2
-\frac{1}{6}\Gamma_6,
\label{donotdrop1}\\
V_7&=&-\frac{1}{24}\Gamma_2^2\Gamma_3+\frac{1}{10}\Gamma_2\Gamma_5
+\frac{1}{12}\Gamma_3\Gamma_4-\frac{1}{7}\Gamma_7,
\label{donotdrop1b}\\
V_8&=&\frac{1}{384}\Gamma_2^4-\frac{1}{32}\Gamma_{2}^2\Gamma_4
-\frac{1}{36}\Gamma_2\Gamma_3^2+\frac{1}{12}\Gamma_2\Gamma_6
+\frac{1}{15}\Gamma_3\Gamma_5+\frac{1}{32}\Gamma_{4}^2
-\frac{1}{8}\Gamma_8,\label{eq77}
\end{eqnarray}
which, in particular, leads to $W_0=W_1=\overline{F}^n$,
\[W_2=\overline{F}^n-\frac{1}{2}\Gamma_2\overline{F}^{n-2}, \qquad
W_3=\overline{F}^n-\frac{1}{2}\Gamma_2\overline{F}^{n-2}
-\frac{1}{3}\Gamma_3\overline{F}^{n-3}.\]
Letting $\ell=1$ and $\alpha=0$, we obtain under the present
assumption that
\begin{equation}\label{eq188} 
\Norm{\prod_{j=1}^nF_j-\overline{F}^n}\leq
\frac{2\ee\,c_1^{}\,\eta_{1}^{}}{1-\sqrt{2\ee\,c_1^{}\,\eta_{1}^{}}},
\quad\mbox{if }\; \eta_1<(2\ee\,c_1)^{-1},
\end{equation}
where 
\begin{equation}\label{eq188b} 
\eta_{1}^{}=\max_{k\,\in\,\set{n}\setminus\set{1}}
\frac{\nu_{k,2}^{}}{k}
\end{equation}
(see comment after Theorem~\ref{thm01})
can be estimated with~(\ref{eq7884}) below.
\item[\textup{(c)}] Let us assume that, for each $j\in\set{n}$, 
$F_j\ll G$ and let $f_j$ be a $G$-density of $F_j$. Set
$\overline{f}=\frac{1}{n}\sum_{j=1}^nf_j$. 
If in Theorem~\ref{thm01}(b) we choose suitable
$B_0,B_1,\dots,B_n\in\{\emptyset,\X\}$, it then follows that, 
for $\ell\in\setn{n}$, 
\begin{eqnarray} 
\eta_{\ell}^{}&\leq&
\frac{1}{4c_1^{}}\min\Big\{2n\int_{\X}(\overline{f}-1)^2\,\dd G,\,
\frac{n^2}{\ell+1}\norm{\overline{F}-G}^2\Big\}\nonumber\\
&&{}+\sum_{j=1}^n\min\Big\{\frac{2}{n}\int_{\X}(f_j-1)^2\,\dd 
G,\,\frac{1}{\ell+1}\norm{F_j-G}^2\Big\}.\label{eq378}
\end{eqnarray}
From the definition of $\eta_\ell^{}$ it is clear that,  if
$G=F_1=\dots =F_n$, then $\eta_\ell^{}=0$ for each $\ell\in\setn{n}$.
The inequalities (\ref{eq378}) and~(\ref{eq512}) reflect this fact. 
Moreover, in view of these bounds, if $G\approx\overline{F}$ in some 
sense and if  the $F_1,\dots,F_n$ are not too different, then a large
$n$  leads to a small bound. Speaking in terms of 
\citet[Introduction]{BHJ92},
our bound contains a magic factor (cf.~Section~\ref{known} above).

\item[\textup{(d)}] If $G=\overline{F}$, then, for each $j\in\set{n}$,
we clearly have $F_j\ll G$ and therefore a $G$-density $f_j$ of 
$F_j$ exists. In this case, (\ref{eq378}) reduces to 
\begin{equation} \label{eq7884}
\eta_{\ell}^{}\leq
\sum_{j=1}^n\min\Big\{\frac{2}{n}\int_{\X}(f_j-1)^2\,\dd 
\overline{F},\,\frac{1}{\ell+1}\norm{F_j-\overline{F}}^2\Big\},
\qquad (\ell\in\set{n}).
\end{equation}
We note that, in (\ref{eq7884}), 
$\int_{\X}(f_j-1)^2\,\dd\overline{F}$ is finite
for all $j\in\set{n}$, which follows from
\[\int_{\X}f_j^2\,\dd\overline{F}=
\int_{\X}f_j\,\dd F_j\leq n\int_{\X}f_j\,\dd \overline{F}=n.\]
\end{itemize} 
\end{remark}

One might ask whether the 
singularity in the right-hand side of~(\ref{eq384b}) can be removed.
The following theorem shows, that this is possible,
if we enlarge the leading absolute constant and replace $\eta_\ell^{}$
with $\eta_0^{}$ (or with $\eta_1^{}$ in the case $G=\overline{F}$).
\begin{theorem}\label{thm02}
Let the notation of Theorem~\ref{thm01} be valid.
\begin{itemize}
\item[\textup{(a)}]
Let $\ell\in\setn{n}$ and let $u_\ell\in(0,\infty)$ be the smallest 
possible constant such that, without any restriction on 
$\eta_0^{}$,
\begin{equation}\label{eq385}
\Norm{\prod_{j=1}^nF_j-W_\ell}\leq u_\ell\,\eta_0^{(\ell+1)/2}.
\end{equation}
We have
\begin{equation}\label{eq3981} 
u_\ell\leq\frac{(2\ee\,c_1^{})^{(\ell+1)/2}}{1-x_\ell},
\end{equation}
where $x_\ell\in(0,1)$ is the unique positive solution of the 
equation $x^{\ell+1}+x/2=1$. By~\textup{(\ref{eq3981})}, we get 
$u_0\leq 5.9$, $u_1\leq 17.3$, $u_2\leq 44.5$, and $u_3\leq107.5$.
\item[\textup{(b)}]
Let $\ell\in\set{n}$  and let $\widetilde{u}_\ell\in(0,\infty)$ be 
the smallest possible constant such that, under the assumption 
$G=\overline{F}$ and without any restriction on $\eta_{1}^{}$, 
\begin{equation}\label{eq227}
\Norm{\prod_{j=1}^nF_j-W_\ell}\leq 
\widetilde{u}_\ell\,\eta_{1}^{(\ell+1)/2}.
\end{equation}
Then we get
\begin{eqnarray}\label{eq199}
\widetilde{u}_1\leq10.94,\quad \widetilde{u}_2\leq31.5,
\quad \widetilde{u}_3\leq82.2, \quad 
\widetilde{u}_\ell\leq \frac{(2\ee\,c_1^{})^{(\ell+1)/2}}{
1-\widetilde{x}_\ell},\quad(\ell\in\set{n}\setminus\set{3}),
\end{eqnarray}
where $\widetilde{x}_\ell\in(0,1)$ is the unique positive solution
of the equation 
$\widetilde{x}^{\ell+1}-\widetilde{x}^2/2+\widetilde{x}=1$.
\end{itemize}
\end{theorem} 
\begin{remark}
\begin{itemize}
\item[\textup{(a)}]
If $\eta_0^{}$, resp.~$\eta_{1}^{}$, is  sufficiently small, the
bounds given in Theorem~\ref{thm02} can be further improved as 
follows from Theorem~\ref{thm01} and Lemma~\ref{la66} below. 
In particular, in the case $G=\overline{F}$,  we have 
(cf.\ proof of Theorem~\ref{thm02})
\begin{eqnarray} 
\Norm{\prod_{j=1}^nF_j-\overline{F}^n}&\leq&
\frac{1}{2}\nu_{2,2}+\Norm{\prod_{j=1}^nF_j-W_2}
\leq (1+\widetilde{u}_2\sqrt{\eta_1^{}})\eta_1^{},\label{eq1771}\\
\Norm{\prod_{j=1}^nF_j-W_2}&\leq&
\frac{1}{3}\nu_{3,3}+\Norm{\prod_{j=1}^nF_j-W_3}\leq 
(\sqrt{3}+\widetilde{u}_3\sqrt{\eta_1^{}})\eta_1^{3/2},\nonumber\\
\Norm{\prod_{j=1}^nF_j-W_3}&\leq&
\frac{1}{8}\nu_{4,2}^2+\Norm{\prod_{j=1}^nF_j-W_4}\leq 
(2+\widetilde{u}_4\sqrt{\eta_1^{}})\eta_1^{2}.\nonumber
\end{eqnarray}
In view of (\ref{eq1771}), one may conjecture that 
$\widetilde{u}_1\geq1$. 
Indeed, this is correct and follows from the simple observation that,
for $\X=\ZZ$, $n\in2\NN$, $F_1=\dots=F_{n/2}=\dirac_0$,
$F_{n/2+1}=\dots=F_n=\dirac_1$, we have
\[\Norm{\prod_{j=1}^nF_j-\overline{F}^n}=
2\Big(1-\binomial{n}{n/2}\frac{1}{2^n}\Big)\longrightarrow2,
\quad(n\to\infty)\]
and, by (\ref{eq227}) and (\ref{eq7884}), 
$\norm{\prod_{j=1}^nF_j-\overline{F}^n}\leq\widetilde{u}_1\,\eta_1
\leq 2\widetilde{u}_1$.
\item[\textup{(b)}] From (\ref{eq227}) and (\ref{eq7884}), 
it follows that
\[\Norm{\prod_{j=1}^nF_j-\overline{F}^n}\leq 2\widetilde{u}_1
\max_{j\in\set{n}}\int_{\X}(f_j-1)^2\,\dd \overline{F}.\]
\item[\textup{(c)}]
It is unclear, whether it is possible to remove the singularity 
in~(\ref{eq384}) for any $\alpha>0$. Indeed, since the denominator 
of the right-hand side of (\ref{eq384}) contains $\eta_l$ and not 
$\eta_{\ell,\alpha}$, we cannot argue as in the proof 
of Theorem~\ref{thm02}.
\end{itemize}
\end{remark}
\begin{example}\label{ex115}
In the situation of Theorem~\ref{thm01}, let us assume that
$F_j=\sum_{r=0}^d p_{j,r}H_r$, ($j\in\set{n}$, $d\in\NN$) and 
$G=\overline{F}=\sum_{r=0}^d\overline{p}_rH_r$,  where
$H_0,\dots,H_d\in\F$, and for $r\in\setn{d}$, 
$p_{j,r}\in[0,1]$ with $\sum_{r=0}^dp_{j,r}=1$
and $\overline{p}_r=\frac{1}{n}\sum_{j=1}^np_{j,r}>0$.  
Then, for each $r\in\setn{d}$, $H_r$ has a 
$\overline{F}$-density $h_r$ and we may assume that
$\sum_{r=0}^d\overline{p}_rh_r=1$. Consequently, $F_j$ has the
$\overline{F}$-density $f_j\defby\sum_{r=0}^dp_{j,r}h_r$, $(j\in\set{n})$.
Using the simple  inequality 
\begin{equation}\label{eq4876} 
\frac{(\sum_{r=0}^da_r)^2}{\sum_{r=0}^da_r'}\leq 
\sum_{r=0}^d\frac{a_r^2}{a_r'},\quad(a_r\in[0,\infty),
a_r'\in(0,\infty)\mbox{ for }r\in\setn{d}),
\end{equation}
we obtain, for $j\in\set{n}$, 
\begin{eqnarray*} 
\int_{\X}(f_j-1)^2\,\dd\overline{F}&=& 
\int_{\X}f_j^2\,\dd\overline{F}-1
=\int_{\X}\frac{(\sum_{r=0}^dp_{j,r}h_r)^2}{
\sum_{r=0}^d\overline{p}_rh_r}\,\dd\overline{F}-1\\
&\leq&\sum_{r=0}^d\frac{p_{j,r}^2}{\overline{p}_r}
\int_{\{h_r>0\}}\frac{h_r^2}{h_r}\,\dd\overline{F}-1
=\sum_{r=0}^d\frac{p_{j,r}^2}{\overline{p}_r}-1
=\sum_{r=0}^d\frac{(\overline{p}_r-p_{j,r})^2}{\overline{p}_r}.
\end{eqnarray*} 
Further, we have 
\[\norm{F_j-\overline{F}}
\leq\sum_{r=0}^d\ab{\overline{p}_r-p_{j,r}}.\]
Therefore, in this context, (\ref{eq7884}) implies that, for
$\ell\in\set{n}$, 
\begin{equation}\label{eq5511}
\eta_{\ell}^{}\leq
\sum_{j=1}^n\min\Big\{2\sum_{r=0}^d\frac{(\overline{p}_r-p_{j,r})^2
}{n\overline{p}_r},\,\frac{1}{\ell+1}\Big(\sum_{r=0}^d
\ab{\overline{p}_r-p_{j,r}}\Big)^2\Big\}.
\end{equation}
Using~(\ref{eq4876}), we get  
\[\frac{(\overline{p}_0-p_{j,0})^2}{\overline{p}_0}
=\frac{(1-\overline{p}_0)(\sum_{r=1}^d(\overline{p}_r-p_{j,r}))^2}{
\overline{p}_0\sum_{r=1}^d\overline{p}_r}
\leq\Big(\frac{1}{\overline{p}_0}-1\Big)
\sum_{r=1}^d\frac{(\overline{p}_r-p_{j,r})^2}{\overline{p}_r}\]
and hence
\begin{equation}\label{eq2644}
\sum_{r=0}^d\frac{(\overline{p}_r-p_{j,r})^2}{n\overline{p}_r}
\leq\sum_{r=1}^d\frac{(\overline{p}_r-p_{j,r})^2}{
n\overline{p}_r\overline{p}_0}.
\end{equation}
We note that (\ref{eq2644}) is non-trivial in the sense that the sum 
on the right-hand side does not contain the summand for $r=0$. 
In view of~(\ref{eq227}), (\ref{eq199}), and~(\ref{eq5511}) 
with $\ell=1$, and~(\ref{eq2644}), 
we see that, in~(\ref{eq3590}), the factor $C_3d$ can be
replaced with $2\widetilde{u}_1$, which in turn is bounded by $21.88$.
We note that, if the $H_r$ are given as in~(\ref{eq1199}) 
then~(\ref{eq5511}) and~(\ref{eq7884}) coincide. 
But if $H_0\approx\dots\approx H_d$ in some sense, 
then~(\ref{eq5511}) can be much worse than~(\ref{eq7884}) and  
should therefore not be used in general.
\end{example}

The next proposition shows that, as claimed above, sometimes
$\eta_{\ell,\alpha}^{}$, $(\alpha>0)$ has a better order 
than~$\eta_{\ell}^{}$. Here, we consider the case of symmetric 
distributions $F_1,\dots,F_n\in\F$ with finite support. 
For simplicity, we assume that $G=\overline{F}$. 
\begin{proposition}\label{prop256}
Let the notation from Theorem~\ref{thm01} hold.
Further, let $b\in\NN$, $x_1,\dots,x_b\in\X\setminus\{0\}$, 
$F_j=p_{j,0}\dirac+\sum_{r=1}^bp_{j,r}(\dirac_{-x_r}+
\dirac_{x_r})\in\F$, $(j\in\set{n})$, and 
$G=\overline{F}
=\overline{p}_{0}\dirac+\sum_{r=1}^b\overline{p}_r
(\dirac_{-x_r}+\dirac_{x_r})$,
where $p_{j,r}\in[0,1]$ with $p_{j,0}+2\sum_{r=1}^b p_{j,r}=1$ and 
$\overline{p}_r=\frac{1}{n}\sum_{j=1}^np_{j,r}>0$, $(r\in\setn{b})$.
For $\ell\in\set{n}$, we then have  
\begin{eqnarray}
\eta_{\ell,1}^{}&\leq&\sum_{j=1}^n\min\Big\{\frac{4}{n^2}
\Big(\frac{(\overline{p}_0-p_{j,0})^2}{2\overline{p}_0^{2}}
+2\sum_{r=1}^b\frac{(\overline{p}_r-p_{j,r})^2}{\overline{p}_r
\overline{p}_0}+\sum_{r=1}^b\frac{(\overline{p}_r-p_{j,r})^2}{
\overline{p}_r^2}\Big),\nonumber\\
&&\hspace{6cm}\frac{1}{(\ell+1)^2}\Big(\ab{\overline{p}_0-p_{j,0}}+
2\sum_{r=1}^b\ab{\overline{p}_r-p_{j,r}}\Big)^2\Big\}.\label{eq181}
\end{eqnarray}
\end{proposition}
We note that, in contrast to~(\ref{eq7884}), the bound 
in~(\ref{eq181}) has the better magic factor $n^{-2}$. Hence, in the
situation of Proposition~\ref{prop256}, estimate (\ref{eq384}) 
with $\alpha=1$ should be preferred over~(\ref{eq384b}).
\section{Numerical examples}\label{numerics}
In what follows, we compare the available bounds in the multinomial 
approximation of the generalized multinomial distribution. 
We assume the notation given in (\ref{eq1199}) with $d=10$.
Further let $\ell=1$.
The following two examples show that the results of the present paper 
can be considerably sharper than the bounds from the literature 
discussed in Section~\ref{known}. 
\begin{example}\label{ex78}
For $j\in\set{n}$, let $p_{j,r}=\binomial{d}{r}q_j^r(1-q_j)^{d-r}$,
$(r\in\setn{d})$ 
be the binomial counting density with number of trials $d$ and 
success probability $q_j=0.4+\frac{1}{(j+9)^a}$, where $a\geq1$. 
Clearly we have $q_j\in(0.4,\,0.5]$ for all $j\in\set{n}$.
We emphasize that, with this definition, $F_j$ is not a binomial 
distribution. Further, if $a$ or $n$ is large, then 
$p_{j,r}$ should be close to $\overline{p}_r$
for a sufficient number of $j\in\set{d}$ and $r\in\setn{d}$, 
so that we expect a small distance 
$\norm{\prod_{j=1}^nF_j-\overline{F}^n}$ here. 
This is reflected in the bounds, given in Table~\ref{tab16}.
\begin{table}[H]\caption{Numerical bounds for the distance in
Example~\ref{ex78}}
\label{tab16}
\begin{tabular*}{\textwidth}{@{\extracolsep{\fill}}ccccccccc}\hline
$n$ & $a$& $C_1$ & $C_2$ & (\ref{loh})  & (\ref{eq1}) & (\ref{eq1b}) & 
(\ref{eq1771})~\&~(\ref{eq5511})   
& (\ref{eq188})~\&~(\ref{eq5511})  \\ \hline
$100$ & $1$  & $111.4$ & $15590.9$& n.a. &  $\geq2$ &  n.a. 
& $0.197438$ & $0.173503$\\
$1000$ & $1$  & $145.7$ & $26444.8$& n.a. &  $\geq2$ &  n.a. 
& $0.026902$ & $0.032981$\\
$100$ & $2$  &  $154.6$ & $29809.2$& n.a. & $0.107737$ & $0.034777$ 
& $0.000366$ & $0.000954$\\
$1000$ & $2$ & $156.3$ & $30455.0$& n.a. & $0.110925$ &  $0.035914$ 
& $0.000037$ & $0.000120$\\ \hline 
\end{tabular*}
\end{table}
\noindent
Note that the bounds for the distance are always rounded up. 
Further, as the distance is always bounded by $2$, larger bounds are
omitted.
The entry ``n.a.'' means ``not available'' and describes a situation, 
where the bound cannot be used, since the respective condition does 
not hold. In all cases, the quantities $C_1$ and $C_2$
(see (\ref{eq6341}) for the definition) are quite 
large, which explains that the condition for (\ref{loh}) is not 
valid here. This is due to the fact that, in each case,
some of the $\overline{p}_r$, $(r\in\setn{d})$ are quite small.
E.g.~see Table~\ref{tab33} for the case $n=100$ and $a=1$.
\begin{table}[H]\caption{Point probabilities of $\overline{F}$ in
Example~\ref{ex78} when $n=100$, $a=1$}\label{tab33}
\begin{tabular*}{\textwidth}{@{\extracolsep{\fill}}c|cccccc}\hline
$r$&$0$&$1$&$2$&$3$&$4$&$5$\\ 
$\overline{p}_r$
&$0.00416$
&$0.03012$
&$0.09851$
&$0.19175$
&$0.24611$
&$0.21781$\\ \hline\hline
$r$&$6$&$7$&$8$&$9$&$10$&\\
$\overline{p}_r$
&$0.13473$
&$0.05757$
&$0.01628$
&$0.00276$
&$0.00021$
& \\ \hline 
\end{tabular*}
\end{table}
\end{example}
In the next example, we discuss a situation, where $(\ref{loh})$ 
gives non-trivial bounds. 
\begin{example}\label{ex2}
For $j\in\set{n}$ and $r\in\setn{d}$, let 
\[p_{j,r}=\frac{1+(j+r)/(b(n+d))}{\sum_{r_1=0}^d(1+(j+r_1)/(b(n+d)))},\]
where $b\geq1$. Similarly as in Example~\ref{ex78}, 
for large $n$ or $b$, we expect good approximation, which indeed 
is reflected in the bounds 
for $\norm{\prod_{j=1}^nF_j-\overline{F}^n}$ given in 
Table~\ref{tab62}.
\begin{table}[H]\caption{Numerical bounds for the distance in
Example~\ref{ex2}}\label{tab62}
\begin{tabular*}{\textwidth}{@{\extracolsep{\fill}}ccccccc}\hline
$n$  & $b$ & (\ref{loh}) & (\ref{eq1}) & (\ref{eq1b}) & 
(\ref{eq1771}) \& (\ref{eq5511})  
& (\ref{eq188}) \& (\ref{eq5511}) \\ \hline
$100$ & $1$ &$0.325253$ & $0.008310$ &  $0.002337$ & $0.000030$ 
& $0.000098$ \\
$1000$ & $1$ &$0.118021$ & $0.000119$ &  $0.000033$ & 
$3.9\times10^{-7}$
& $1.5\times10^{-6}$ \\
$100$ & $2$ 
& $0.112763$ & $0.000978$ & $0.000267$ & $3.3\times10^{-6}$ 
& $1.2\times10^{-5}$\\
$1000$ & $2$ 
& $0.040581$ & $0.000014$ & $3.8\times10^{-6}$ 
& $4.4\times10^{-8}$ & $1.7\times10^{-7}$\\
\hline
\end{tabular*}
\end{table}
\noindent
In contrast to Example~\ref{ex78}, in each case the
values $\overline{p}_r$, $(r\in\setn{d})$ are quite similar, 
which implies that the condition for (\ref{loh}) is valid.
E.g.~see Table~\ref{tab19} for the case $n=100$ and $b=1$.
\begin{table}[H]\caption{Point probabilities of $\overline{F}$ in 
Example~\ref{ex2} when $n=100$, $b=1$}\label{tab19}
\begin{tabular*}{\textwidth}{@{\extracolsep{\fill}}c|cccccc}\hline
$r$&$0$&$1$&$2$&$3$&$4$&$5$\\ 
$\overline{p}_r$
&$0.08807$
&$0.08864$
&$0.08921$
&$0.08978$
&$0.09034$
&$0.09091$\\ \hline\hline
$r$&$6$&$7$&$8$&$9$&$10$&\\
$\overline{p}_r$
&$0.09148$
&$0.09204$
&$0.09261$
&$0.09318$
&$0.09374$&
\\\hline
\end{tabular*}
\end{table}
\end{example}
In what follows, we discuss an example, where 
the distance can actually be evaluated. 
\begin{example}
Suppose now that, in Example~\ref{ex78}, we change the measures $H_r$ 
to $H_r=\dirac_r$ on $\RR$ for $r\in\setn{d}$, i.e.\ all
distributions $F_1,\dots,F_n,\overline{F}$ are one-dimensional. 
Then, using a computer, it is not difficult to get the exact numerical
value for the distance, see Table~\ref{tab52}.
\begin{table}[H]\caption{Exact numerical values for the distance
(cf.\ with Table~\ref{tab16})}\label{tab52}
\begin{tabular*}{\textwidth}{@{\extracolsep{\fill}}c|cccc}\hline
$n$ & $100$ & $1000$ & $100$ & $1000$\\
$a$ &  $1$  & $1$    & $2$   & $2$\\
$\norm{\prod_{j=1}^nF_j-\overline{F}^n}$ 
& $0.007152$
& $0.001653$
&$5.9\times10^{-5}$
&$7.6\times10^{-6}$
\\ \hline
\end{tabular*}
\end{table}
\noindent
A basic property of the total variation distance tells us that,
for distributions $\widetilde{H}_r\in\F$, $(r\in\setn{d})$
in the case of a general measurable Abelian group, we have 
\begin{equation}\label{eq723}
\Norm{\prod_{j=1}^n\Big(\sum_{r=0}^dp_{j,r}\widetilde{H}_r\Big)
-\Big(\sum_{r=0}^d\overline{p}_{r}\widetilde{H}_r\Big)^n}
\leq \Norm{\prod_{j=1}^n\Big(
p_{j,0}\dirac+\sum_{r=1}^dp_{j,r}\dirac_{\unitvec{r}}\Big)
-\Big(\overline{p}_{0}\dirac
+\sum_{r=1}^d\overline{p}_{r}\dirac_{\unitvec{r}}\Big)^n}.
\end{equation}
This can easily be seen by writing the difference of the measures 
on the left-hand side as a polynomial in 
$\widetilde{H_r}$, $(r\in\setn{d})$ and then applying the
triangle inequality. As a consequence of (\ref{eq723}), 
each bound from Table~\ref{tab16} is valid here as well.
A comparison shows that the bounds are getting closer to the 
actual distance as $n$ or $a$ is becoming large. For example, 
the bounds from (\ref{eq1771}) \& (\ref{eq5511}) are about 
$27.6$, $16.3$, $6.2$, and $4.9$ times higher, respectively, 
than the values from Table~\ref{tab52}.

We can apply this idea to Example~\ref{ex2} as well:
if we again change the measures $H_r$ to $H_r=\dirac_r$ for 
$r\in\setn{d}$, we get the exact values of Table~\ref{tab69}.
A comparison with Table~\ref{tab62} shows that 
the bounds from (\ref{eq1771}) \& (\ref{eq5511}) are about 
$4.8$ to $4.0$ times higher than these values.
\begin{table}[H]\caption{Exact numerical values for the distance
(cf.\ with Table~\ref{tab62})}\label{tab69}
\begin{tabular*}{\textwidth}{@{\extracolsep{\fill}}c|cccc}\hline
$n$ & $100$ & $1000$ & $100$ & $1000$\\
$b$ &  $1$  & $1$    & $2$   & $2$\\
$\norm{\prod_{j=1}^nF_j-\overline{F}^n}$ 
& $6.3\times10^{-6}$
& $9.1\times10^{-8}$
&$7.4\times10^{-7}$
&$1.1\times10^{-8}$
\\ \hline
\end{tabular*}
\end{table}
\end{example}
\section{Proofs}\label{proofs}
\subsection{Smoothness estimates for convolutions}\label{aux}
In what follows, we use the standard multi-index  notation: For
$z=(z_1,\dots,z_d)\in\CC^d$, $(d\in\NN)$ and
$w=(w_1,\dots,w_d)\in\Zpl^d$, we set  $z^w=\prod_{r=1}^d z_r^{w_r}$, 
$\vecsum{w}=\sum_{r=1}^dw_r$, and $w!=\prod_{r=1}^d w_r!$. Similarly,
for  $V=(V_1,\dots,V_d)\in\M^d$, set $V^w=\prod_{r=1}^dV_r^{w_r}$. For
$v,w\in\Zpl^d$, we write $v\leq w$ in the case that  $v_r\leq w_r$ for
all $r\in\set{d}$; let  $v\wedge w=(v_1\wedge w_1,\dots,v_d\wedge
w_d)$. Sums over $v$, $\widetilde{v}$, and $w$ are taken over subsets
of  $\Zpl^d$ as indicated. The following lemma is a counterpart of
Lemma~5 in \citet{roo01}.
\begin{lemma}\label{la154}
Let $k,n\in\Zpl$, $d\in\NN$, and $a_v\in\RR$  for $v\in\Zpl^d$ with $\ab{v}=k$. Let $X=(X_r)_{r\in\set{d}}$ be a  random vector in $\RR^d$ with  
$\Expect[(\sum_{r=1}^d\ab{X_r})^k]<\infty$ and put 
$X_0=\sum_{r=1}^dX_r$. Let  $p=(p_r)_{r\in\set{d}}\in(0,1)^d$ such
that $p_0=1-\sum_{r=1}^d p_r\in(0,1)$. Further, let
$H=(H_r)_{r\in\set{d}}\in\F^d$, $H_0\in\F$,
\[G=\sum_{r=0}^d p_r H_r\in\F,\quad
U_1=\sum_{\vecsum{v}=k}\frac{a_v}{v!}\prod_{r=1}^d(H_r-H_0)^{v_r},
\quad
U_2=\Expect\Big(\sum_{r=1}^d X_r(H_r-H_0)\Big)^k,\]
where, in the definition of $U_2$, the expectation is defined 
setwise. Then we have
\begin{eqnarray} 
\norm{U_1\,G^n}&\leq&\frac{\sqrt{n!}}{\sqrt{(n+k)!}}
\Big(\sum_{\vecsum{w}\leq k}\frac{w!(k-\vecsum{w})!}{p^w\,
p_0^{k-\vecsum{w}}}
\Big[\sum_{\vecsum{v}=k}\frac{a_{v}}{v!}\prod_{r=1}^d
\binomial{v_r}{w_r}\Big]^2\Big)^{1/2},\label{eq551a}\\
\norm{U_2\,G^n}&\leq&\binomial{n+k}{k}^{-1/2}\Big(\Expect\Big(\sum_{
r=0}^d\frac{X_r Y_r}{p_r}\Big)^k\Big)^{1/2},\label{eq551}
\end{eqnarray}
where the random vector $Y=(Y_r)_{r\in\set{d}}$ is an independent 
copy of~$X$ and $Y_0=\sum_{r=1}^dY_r$.
\end{lemma}
\Proof
Let 
\[\muc(w,n,p)=\left\{\begin{array}{ll}
\frac{n!}{w!\,(n-\vecsum{w})!}\,p^w\,p_0^{n-\vecsum{w}},&
\quad\mbox{if } w\in\Zpl^d,\;\vecsum{w}\leq n,\\
0,&\quad \mbox{otherwise}
\end{array}\right.\]
denote the multinomial counting density with parameters~$n$  and~$p$.
For $f\,:\,\ZZ^d\longrightarrow\RR$ and $r\in\set{d}$,  let $\Delta_r
f\,:\,\ZZ^d\longrightarrow\RR$ with  $(\Delta_r f)(w)=f(w-
\unitvec{r})-f(w)$
for $w\in\ZZ^d$.  Products and powers of $\Delta$-operators are
understood in the  sense of composition. Further, let $\Delta_r^0f=f$.
Clearly, 
$\Delta_{r_1}\Delta_{r_2}f=
\Delta_{r_2}\Delta_{r_1}f$ for $r_1,r_2\in
\set{d}$. For $v\in\Zpl^d$, let
$\Delta^vf=\Delta_1^{v_1}\cdots\Delta_d^{v_d}f$. 
We set $\Delta^v\muc(w,n,p)=(\Delta^v\muc(\cdot,n,p))(w)$
for $w\in\Zpl^d$. 
We use the following properties of the multinomial 
distribution (see \citet[formulas (20), (21), and (4)]{roo01}):
For $v\in\Zpl^d$,
\begin{equation}\label{e3}
\sum_{\vecsum{w}\leq n+\vecsum{v}}\Delta^v\muc(w,n,p)\,
H^{w}\,H_0^{n+\vecsum{v}-\vecsum{w}}=G^n\prod_{r=1}^d(H_r-H_0)^{v_r}
\end{equation}
and, for $v,w\in\Zpl^d$, 
\begin{equation}\label{eq225}
\Delta^v\muc(w,n,p)=\Kr(v;w,n+\vecsum{v},p)\muc(w,n+\vecsum{v},p)
\frac{v!\,n!}{(n+\vecsum{v})!\,p^v\,p_0^{\vecsum{v}}},
\end{equation}
where 
\begin{equation}\label{kra1}
\Kr(v;w,n,p)=
\sum_{\widetilde{v}\leq v}
{\binomial{n-\vecsum{w}}{\vecsum{v-\widetilde{v}}}}\frac{
\vecsum{v-\widetilde{v}}!\,(-p)^{v-\widetilde{v}}
\,p_0^{\vecsum{\widetilde{v}}}}{(v-\widetilde{v})!}
\prod_{r=1}^d{\binomial{w_r}{\widetilde{v}_r}}
\end{equation}
is a Krawtchouk polynomial of degree $v$. Note that there is another
set of Krawtchouk polynomials, which forms, together with the one
from~(\ref{kra1}), a bi-orthogonal system of polynomials with respect
to the multinomial distribution (see also \citet{tra89}). From the
more general Lemma~2 in \citet{roo01}, it follows that,  for
$v,\widetilde{v}\in\Zpl^d$ with $\vecsum{v}=\vecsum{\widetilde{v}}$,
we have
\begin{eqnarray}
\lefteqn{\sum_{\vecsum{w}\leq n+\vecsum{v}}\muc(w,n+\vecsum{v},p)\,
\Kr(v;w,n+\vecsum{v},p)\,
\Kr(\widetilde{v};w,n+\vecsum{v},p)}\nonumber\\
&\hspace{6cm}=&\sum_{w\leq v\wedge \widetilde{v}}
\frac{(n+\vecsum{v})!\,\vecsum{v-w}!\,p^{v+\widetilde{v}-w}\,
p_0^{\vecsum{w+v}}}{w!\,n!\,(v-w)!\,(\widetilde{v}-w)!}.\label{eq332}
\end{eqnarray}
We note that the right-hand side of~(\ref{eq332}) is always  positive,
which shows that, if $d\geq2$, then  the Krawtchouk polynomials given
above are not orthogonal with  respect to the multinomial
distribution. However, we do not need such a property. 
Using~(\ref{e3}), (\ref{eq225}), Cauchy's inequality, we now obtain
\begin{eqnarray*}
\norm{U_1\,G^n}
&=&\Norm{\sum_{\vecsum{v}=k}\frac{a_v}{v!}\,
\sum_{\vecsum{w}\leq n+\vecsum{v}}\Delta^{v}\muc(w,n,p)\,H^{w}\,
H_0^{n+\vecsum{v}-\vecsum{w}}}\\
&\leq&\frac{n!}{(n+k)!}\sum_{w\in\Zpl^d}\muc(w,n+k,p)
\Ab{\sum_{\vecsum{v}=k}\frac{a_v}{
p^v\,p_0^k}\Kr(v;w,n+k,p)}\\
&\leq&\frac{n!}{(n+k)!}\Big(\sum_{w\in\Zpl^d}
\muc(w,n+k,p)\Big[\sum_{\vecsum{v}=k}
\frac{a_v}{p^v\,p_0^k}\Kr(v;w,n+k,p)\Big]^2
\Big)^{1/2}\bydef T.
\end{eqnarray*}
Using~(\ref{eq332}), we get 
\begin{eqnarray*}
T&=&\frac{n!}{(n+k)!}\Big(\sum_{\vecsum{v}=k}
\sum_{\vecsum{\widetilde{v}}=k}
\frac{a_v\,a_{\widetilde{v}}}{p^{v+\widetilde{v}}\,
p_0^{2k}}\sum_{w\leq v\wedge \widetilde{v}}\frac{(n+k)!\,
(k-\vecsum{w})!\,
p^{v+\widetilde{v}-w}\,p_0^{\vecsum{w}+k}}{w!\,n!\,(v-w)!\,
(\widetilde{v}-w)!}\Big)^{1/2}\\
&=&\frac{\sqrt{n!}}{\sqrt{(n+k)!}}\Big(\sum_{\vecsum{w}\leq k}
\frac{w!(k-\vecsum{w})!}{p^w\,p_0^{k-\vecsum{w}}}
\Big[\sum_{\vecsum{v}=k}\frac{a_v}{v!}
\prod_{r=1}^d\binomial{v_r}{w_r}
\sum_{\vecsum{\widetilde{v}}=k}\frac{a_{\widetilde{v}}}{
\widetilde{v}!}
\prod_{r=1}^d\binomial{\widetilde{v}_r}{w_r}\Big]\Big)^{1/2}\\
&=&\frac{\sqrt{n!}}{\sqrt{(n+k)!}}\Big(\sum_{\vecsum{w}\leq k}
\frac{w!(k-\vecsum{w})!}{p^w\,p_0^{k-\vecsum{w}}}\Big[\sum_{
\vecsum{v}=k}
\frac{a_{v}}{v!}\prod_{r=1}^d\binomial{v_r}{w_r}\Big]^2\Big)^{1/2}.
\end{eqnarray*}
Inequality~(\ref{eq551a}) is shown. 
Since $U_2=\sum_{\vecsum{v}=k}\frac{k!}{v!}\,\Expect\big[X^v\big]
\prod_{r=1}^d(H_r-H_0)^{v_r}$,~(\ref{eq551a}) gives
\begin{eqnarray*}
\norm{U_2\,G^n}&\leq&
\frac{\sqrt{n!}\,k!}{\sqrt{(n+k)!}}\Big(\sum_{\vecsum{w}\leq k}
\frac{w!(k-\vecsum{w})!}{p^w\,p_0^{k-\vecsum{w}}}
\Big[\sum_{\vecsum{v}=k}
\frac{\Expect{X^v}}{v!}\prod_{r=1}^d\binomial{v_r}{w_r}\Big]^2
\Big)^{1/2}\\
&=&\frac{\sqrt{n!}\,k!}{\sqrt{(n+k)!}}\Big(\sum_{\vecsum{w}\leq k}
\frac{w!(k-\vecsum{w})!}{p^w\,p_0^{k-\vecsum{w}}}
\Expect\Big[\sum_{\vecsum{v}=k}\frac{X^{v}}{v!}
\prod_{r=1}^d\binomial{v_r}{w_r}
\sum_{\vecsum{\widetilde{v}}=k}\frac{Y^{\widetilde{v}}}{
\widetilde{v}!}\prod_{r=1}^d\binomial{\widetilde{v}_r}{w_r}
\Big]\Big)^{1/2}.
\end{eqnarray*}
For $\vecsum{w}\leq k$, we have
\[\sum_{\vecsum{v}=k}\frac{X^v}{v!}\prod_{r=1}^d\binomial{v_r}{
w_r}=\frac{X^w\,X_0^{k-\vecsum{w}}}{w!\,(k-\vecsum{w})!}.\]
Indeed, this follows from the identity theorem for power series 
taking into account the following equality of the corresponding 
generating functions
\[\sum_{k=0}^\infty\Big[\sum_{\vecsum{v}=k}\frac{X^{v}}{v!}
\prod_{r=1}^d\binomial{v_r}{w_r}\Big]z^k
=\frac{z^{\vecsum{w}}X^w}{w!}\ee^{zX_0}=\sum_{k=\vecsum{w}}^\infty
\Big[\frac{X^w\,X_0^{k-\vecsum{w}}}{w!\,(k-\vecsum{w})!}\Big]z^k,
\qquad (z\in\CC).\]
From the above, we get
\begin{eqnarray*}
\norm{U_2\,G^n}
&\leq&\binomial{n+k}{k}^{-1/2}\Big(\Expect\Big[\sum_{\vecsum{w}\leq k}
\frac{k!}{w!\,(k-\vecsum{w})!}\Big(\frac{X_0Y_0}{p_0}\Big)^{
k-\vecsum{w}}\prod_{r=1}^d\Big(\frac{X_r Y_r}{p_r}\Big)^{
w_r}\Big]\Big)^{1/2}\\
&=&\binomial{n+k}{k}^{-1/2}\Big(\Expect\Big(\sum_{r=0}^d\frac{
X_r Y_r}{p_r}\Big)^k\Big)^{1/2},
\end{eqnarray*}
which completes the proof of~(\ref{eq551}).~\hfill\qed\bigskip

The following lemma is an important application of Lemma~\ref{la154}
and generalizes formula~(37) in \citet{roo00}.
Another application is given in the proof of 
Proposition~\ref{prop256}, see Section~\ref{remproofs} below.
\begin{lemma} \label{la239}
Let $k\in\NN$, $n\in\Zpl$, $G\in\F$, and $U\in\M$, where we assume 
that $\tvm{U}\ll G$ and that $U(\X)=0$;
let $f^{\pm}$ denote any  Radon-Nikodym densities of~$U^\pm$
with respect to $G$ and put $f=f^+-f^-$. Then 
\begin{equation}\label{eq18798}
\norm{U^k\,G^n}\leq \binomial{n+k}{k}^{-1/2}\Big(\int
f^2\,\dd G\Big)^{k/2}.
\end{equation}
\end{lemma}
\Proof If $\int f^2\,\dd G=\infty$, then~(\ref{eq18798}) is
trivial. In what follows, we assume that $\int f^2\,\dd G<\infty$.
Let $\varepsilon\in(0,1)$ be fixed. Then 
$a_{j,\varepsilon}^\pm\in[0,\infty)$, 
$(j\in\Zpl)$ and pairwise disjoint $B_{j,\varepsilon}\in\A$, 
$(j\in\Zpl)$ exist such that 
\[\bigcup_{j=0}^\infty B_{j,\varepsilon}=\X,\quad  
f_{\varepsilon}^\pm\defby\sum_{j=0}^\infty
a_{j,\varepsilon}^\pm\indicator(B_{j,\varepsilon}),
\quad\mbox{and}\quad
0\leq f^\pm-f_{\varepsilon}^\pm\leq\varepsilon.\]
Here $\indicator(A)$ is the indicator function of a set $A$. 
Let $U_{\varepsilon}^\pm$ be the measures on $(\X,\A)$ with 
$G$-densities $f_{\varepsilon}^\pm$. This implies that 
$U_{\varepsilon}^\pm=\sum_{j=0}^\infty q_{j,\varepsilon}^\pm
\,H_{j,\varepsilon}$,
where, for $j\in\Zpl$, 
\[q_{j,\varepsilon}^\pm=a_{j,\varepsilon}^\pm\,G(B_{j,\varepsilon})
\quad
\mbox {and }\quad
H_{j,\varepsilon}=\left\{
\begin{array}{ll} 
G(B_{j,\varepsilon}\cap \,\cdot\,)/G(B_{j,\varepsilon}),
&\mbox{ if } G(B_{j,\varepsilon})>0, \\
\dirac,&\mbox{ otherwise.} 
\end{array}\right.\]
Set 
$q_{j,\varepsilon}=q_{j,\varepsilon}^+-q_{j,\varepsilon}^-$,
$f_{\varepsilon}=f_{\varepsilon}^+-f_{\varepsilon}^-$, and
$U_{\varepsilon}=U_{\varepsilon}^+-U_{\varepsilon}^-$. 
We note that the latter equality indeed indicates the Hahn-Jordan
decomposition of $U_{\varepsilon}$. Then
$\norm{U_{\varepsilon}}\leq\norm{U}$ and 
$\norm{U-U_{\varepsilon}}=\int\big((f^+-f_{\varepsilon}^+)
+(f^--f_{\varepsilon}^-)\big)\,\dd G\leq\varepsilon$, giving 
\[\norm{U^k-U_{\varepsilon}^k}\leq k\,\norm{U}^{k-1}
\norm{U-U_{\varepsilon}}\leq k\,\norm{U}^{k-1}\varepsilon.\]
Hence 
\[\norm{U^k\,G^n}\leq\norm{(U^k-U_{\varepsilon}^k)G^n}
+\norm{U_{\varepsilon}^k\,G^n}
\leq k\,\norm{U}^{k-1}\varepsilon+\norm{U_{\varepsilon}^kG^n}.\]
Since $\int f\,\dd G=U(\X)=0$, we have 
$\Ab{\sum_{j=0}^\infty q_{j,\varepsilon}}
=\Ab{\int(f_\varepsilon^+-f_\varepsilon^-)\,\dd G}\leq
\norm{U-U_\varepsilon}\leq\varepsilon$, and therefore, for each
$m\in\NN$,
\[\Norm{\Big(\sum_{j=0}^\infty q_{j,\varepsilon}H_{j,\varepsilon}
\Big)^k-\Big(\sum_{j=0}^m q_{j,\varepsilon}(H_{j,\varepsilon}-
H_{0,\varepsilon})\Big)^k}\leq 
k\Big(\varepsilon+2\sum_{j=m+1}^\infty\ab{q_{j,\varepsilon}}\Big)
\Big(2\sum_{j=0}^\infty \ab{q_{j,\varepsilon}}\Big)^{k-1}.\]
Hence,  we obtain
{\arraycolsep2pt
\begin{eqnarray*}
\norm{U_\varepsilon^k\,G^n}
&=&\Norm{
\Big(\sum_{j=0}^\infty q_{j,\varepsilon}H_{j,\varepsilon}\Big)^k
\Big(\sum_{j=0}^\infty G(B_{j,\varepsilon})H_{j,\varepsilon}\Big)^n}\\
&\leq&k\Big(\varepsilon+2\sum_{j=m+1}^\infty\ab{q_{j,\varepsilon}}\Big)
\Big(2\sum_{j=0}^\infty \ab{q_{j,\varepsilon}}\Big)^{k-1}
+\Norm{\Big(\sum_{j=1}^m q_{j,\varepsilon}(H_{j,\varepsilon}-
H_{0,\varepsilon})\Big)^k
\Big(\sum_{j=0}^\infty G(B_{j,\varepsilon})H_{j,\varepsilon}\Big)^n}.
\end{eqnarray*}}%
From (\ref{eq551}), it follows that
the norm term on the right-hand side is bounded from above by
\begin{eqnarray*}
\binomial{n+k}{k}^{-1/2}
\Big(\sum_{j\in\setn{m}:\,G(B_{j,\varepsilon})>0}
\frac{q_{j,\varepsilon}^2}{G(B_{j,\varepsilon})}\Big)^{k/2}
&=&\binomial{n+k}{k}^{-1/2}\Big(\sum_{j=0}^m
(a_{j,\varepsilon}^+-a_{j,\varepsilon}^-)^2
G(B_{j,\varepsilon})\Big)^{k/2}\\
&\leq&\binomial{n+k}{k}^{-1/2}
\Big(\int f_\varepsilon^2\,\dd G\Big)^{k/2}.
\end{eqnarray*}
Letting $m\to\infty$, we obtain 
\[\norm{U_\varepsilon^k\,G^n}\leq k\varepsilon\Big(2\sum_{j=0}^\infty
\ab{q_{j,\varepsilon}}\Big)^{k-1}+\binomial{n+k}{k}^{-1/2}
\Big(\int f_\varepsilon^2\,\dd G\Big)^{k/2}.\]
Since 
\[\Ab{\int (f^2-f_\varepsilon^2)\,\dd G}
\leq\int\ab{f-f_\varepsilon}(\ab{f}+\ab{f_\varepsilon})\,\dd G \leq
\varepsilon(\norm{U}+\norm{U_\varepsilon})
\leq2\varepsilon\norm{U},\]
we obtain (\ref{eq18798}) by letting $\varepsilon\to0$. 
This completes the proof.~\hfill\qed\bigskip

It may happen that the assumption in Lemma~\ref{la239} does not hold
directly. However, this can sometimes be overcome by shifting $U$. 
The following corollary is needed in the proof of Theorem~\ref{thm01}.
\begin{corollary} \label{cor227}
Let $n\in\Zpl$, $G\in\F$, $U_1,U_2\in\M$, and $U=U_1+U_2$.
We assume that $\tvm{U_2}\ll G$ and that both $U_2^\pm\neq0$. Put
$\widetilde{U}_2^\pm=U_2^\pm/\norm{U_2^\pm}$. Let 
$f^{\pm}$ denote any Radon-Nikodym densities of~$\widetilde{U}_2^\pm$
with respect to $G$ and set $f=f^+-f^-$. Then 
\begin{equation}\label{eq1879}
\norm{U\,G^n}\leq \norm{U_1}+\ab{U_2(\X)}+\frac{\norm{U_2^+}\wedge
\norm{U_2^-}}{\sqrt{n+1}}\Big(\int f^2\,\dd G\Big)^{1/2}.
\end{equation}
\end{corollary}
\Proof The assertion easily follows from the triangle inequality,
Lemma~\ref{la239},  and the simple fact that 
$U_2=U_2(\X)\widetilde{U}^{\tau}_2
+(\norm{U_2^+}\wedge\norm{U_2^-})
(\widetilde{U}_2^+-\widetilde{U}_2^-)$, where $\tau$ denotes $+$ or 
$-$ according to whether $\norm{U_2^+}>\norm{U_2^-}$ or 
not.~\hfill\qed
\begin{remark}\label{rem62}
\begin{itemize}
\item[(a)]
Let the assumptions of Corollary~\ref{cor227} hold. If $\mu$ is a
$\sigma$-finite measure on $\X$ and if 
$G\ll\mu$, then $G$ and $\widetilde{U}_2^\pm$ have $\mu$-densities 
$v$ and $g^\pm$, say, and, letting $g=g^+-g^-$, we can write 
$\int f^2\,\dd G=\int_{\{v>0\}}g^2v^{-1}\,\dd\mu$.
\item[(b)] 
Sometimes it is useful to simplify further the bound~(\ref{eq1879}) 
by using the following inequality:
$(\norm{U_2^+}\wedge\norm{U_2^-})^2\int f^2\,\dd G
\leq\int h^2\,\dd G$, where $h=h^+-h^-$ and $h^\pm$ denote any
$G$-densities of~$U_2^\pm$. Indeed, this follows from the
representation 
\[\int f^2\,\dd G=\int_{A}\frac{(h^+)^2}{\norm{U_2^+}^2}\,\dd G
+\int_{A^\cc} \frac{(h^-)^2}{\norm{U_2^-}^2}\,\dd G,\]
whenever $A\in\A$ with $U_2^-(A)=U_2^+(A^\cc)=0$.
\end{itemize}
\end{remark}

The next corollary is an extension of Lemma~\ref{la239} to compound 
distributions and may be particularly useful in the compound Poisson
approximation. 
\begin{corollary}\label{cor25}
Let $k\in\NN$, $G\in\F$, and $U\in\M$, where we assume  that
$\tvm{U}\ll G$ and that $U(\X)=0$; let $f^{\pm}$ denote any
Radon-Nikodym densities of~$U^\pm$ with respect to $G$ 
and put $f=f^+-f^-$. Let $N$ be a random
variable in $\Zpl$ and $\varphi(z)=\Expect[z^N]$, $(z\in\CC,\,
\ab{z}\leq1)$ be its generating function. 
Set $\varphi(G)=\Expect[G^N]\in\F$, where the expectation is 
defined setwise. Then we have
\begin{equation}\label{eq18497}
\norm{U^k\varphi(G)}\leq 
\Big(k\int_0^1x^{k-1}\varphi(1-x)\,\dd x\Big)^{1/2}
\Big(\int f^2\,\dd G\Big)^{k/2}.
\end{equation}
If N has Poisson distribution 
$\expo{t(\dirac_1-\dirac)}$ with $t\in(0,\infty)$, then
\begin{equation}\label{eq1843}
\norm{U^k\varphi(G)}\leq \frac{1}{t^{k/2}}\,\sqrt{k!\,\Prob(N\geq k)}
\Big(\int f^2\,\dd G\Big)^{k/2}.
\end{equation}
\end{corollary}
\Proof 
Using the triangle inequality, Lemma~\ref{la239}, and Jensen's 
inequality, we obtain 
\[\norm{U^k\varphi(G)}\leq\Expect\norm{U^kG^N}
\leq\Big(\Expect\binomial{N+k}{k}^{-1}\Big)^{1/2}
\Big(\int f^2\,\dd G\Big)^{k/2}.\]
The integral representation of the beta function implies that
$\Expect\binomial{N+k}{k}^{-1}=k\int_0^1x^{k-1}\varphi(1-x)\,\dd x$,
which, in turn, leads to (\ref{eq18497}). Inequality (\ref{eq1843}) 
easily follows from (\ref{eq18497}) and the series representation 
of the lower incomplete gamma function.~\hfill\qed

We note that (\ref{eq1843}) is comparable to previous results of 
\citet[Lemma~2]{roo03} but is however much better because of the more 
general assumptions used in Corollary~\ref{cor25}.
\subsection{A general lemma}\label{genla}
The results of Section~\ref{mainres} are based on the following 
general lemma. Here, a distribution $G\in\F$ is called infinitely 
divisible if, for each $n\in\NN$, there exists a $G_n\in\F$ such that 
$G_n^n=G$. We note that, in general, such a $n$-th root $G_n$ 
need not be unique (see \citet[proof of Theorem~3.5.15, 
pp.~222--223]{hey77}); let $G^{1/n}$ denote any fixed $n$-th root 
of~$G$.
\begin{lemma} \label{la1a}
Let $n\in\NN$, $F_1,\dots,F_n,G\in\F$, $L_1,\dots,L_n\in\M$. Set 
$\overline{L}=\frac{1}{n}\sum_{j=1}^nL_j$,
$K_j=F_j\ee^{-L_j}$, $(j\in\set{n})$, $K_0=G\ee^{-\overline{L}}$,
\begin{eqnarray*} 
&\displaystyle V_k=\sum_{J\subseteq\,\set{n}:\;\card{J}=k}
\prod_{j\in J}(K_j-K_0),\quad (k\in\setn{n}),\quad
\quad W_{\ell}=\sum_{k=0}^\ell V_k\,G^{n-k}\ee^{k\overline{L}},
\quad (\ell\in\setn{n}),&\\
&\displaystyle 
M_{j,k}=\left\{\begin{array}{ll}
(K_j-K_0)(G^{n-k})^{1/k}\ee^{\overline{L}},&
\quad\mbox{if $G$ is infinitely divisible},\\
(K_j-K_0)G^{\floor{(n-k)/k}}\ee^{\overline{L}},&
\quad\mbox{otherwise,} 
\end{array}\right.
\quad (j,k\in\set{n}),&\\
&\displaystyle \nu_{k,m}^{}=\sum_{j=1}^n\norm{M_{j,k}}^m,\quad
\widetilde{\nu}_k^{}=\Norm{\sum_{j=1}^nM_{j,k}},
\quad (k\in\set{n},m\in\NN).&
\end{eqnarray*}
Let $c_1^{}$ be defined as in Theorem~\ref{thm01}.
Then, for $\ell\in\setn{n}$,
\[\Norm{\prod_{j=1}^nF_j-W_\ell}\leq
\sum_{k=\ell+1}^n
\Big[\Big(\frac{2\ee\,c_1^{}}{k}
\Big(\frac{\widetilde{\nu}_k^2}{4c_1^{}}+\nu_{k,2}^{}\Big)
\Big)^{k/2}\wedge\frac{\nu_{k,1}^k}{k!}\Big].\]
\end{lemma}

The following two examples show possible applications of 
Lemma~\ref{la1a}. As a byproduct, results in the compound Poisson 
approximations can be derived.
\begin{example}
If we consider the case $L_1=\dots=L_n=0$, we see that 
Theorem~\ref{thm01}(a) is a direct consequence of Lemma~\ref{la1a},
(cf.~proof of Theorem~\ref{thm01}).
\end{example}
\begin{example}\label{ex021} 
Suppose that, for $j\in\set{n}$, $H_j\in\F$, $p_j\in[0,1]$,
$L_j=p_j(H_j-\dirac)$, and $F_j=\dirac+L_j$. Put
$\overline{L}=n^{-1}\sum_{j=1}^nL_j$ and $G=\ee^{\overline{L}}$. 
Then Lemma~\ref{la1a} implies that
\begin{equation}\label{eq27} 
\Norm{\prod_{j=1}^nF_j-G^n}\leq\sum_{k=1}^n\frac{1}{k!}
\Big(\sum_{j=1}^n\norm{M_{j,k}}\Big)^k,
\end{equation}
where, for $j,k\in\set{n}$,
\begin{eqnarray*} 
K_j&=&(\dirac+p_j(H_j-\dirac))\ee^{-p_j(H_j-\dirac)},\qquad 
K_0=\dirac,\\
M_{j,k}&=&(K_j-K_0)(G^{n-k})^{1/k}\ee^{\overline{L}}
=((\dirac+L_j)\ee^{-L_j}-\dirac)\Expo{\frac{n}{k}\overline{L}}.
\end{eqnarray*}
In principle, (\ref{eq27}) is the same as estimate (26) in 
\citet{roo03}. The approach used there is based on a slight 
modification of an expansion due to \citet{ker64}. It is however not
sufficient  to get the results of the present paper.
\end{example}

For the proof of Lemma~\ref{la1a}, we use formal power series 
over~$\M$. In the following lemma, some basic properties in 
connection with the norm on $\M$ are summarized. The proof is simple 
and therefore omitted.
\begin{lemma} \label{la176}
For $n\in\NN$ and $k\in\set{n}$, let 
$\psi_k^{(0)}(z)=\sum_{j=0}^\infty W_{j,k}\,z^j$,
$(W_{j,k}\in\M)$ be a formal power series over $\M$ with variable $z$
and let $\coeff{z^j}{\psi_k^{(0)}(z)}$ be its $j$th coefficient 
$W_{j,k}$. Further, consider the formal power series 
$\psi_k^{(1)}(z)=\sum_{j=0}^\infty \norm{W_{j,k}}\,z^j$
and $\psi_k^{(2)}(z)=\sum_{j=0}^\infty a_{j,k}\,z^j$ for
$a_{j,k}\in[\norm{W_{j,k}},\infty)$ and $k\in\set{n}$.
Then, for $j\in\Zpl$, 
\begin{eqnarray*} 
&\displaystyle \norm{\coeff{z^j}{\psi_1^{(0)}(z)}}
=\coeff{z^j}{\psi_1^{(1)}(z)},&\\
&\displaystyle \Norm{\Coeff{z^j}{\prod_{k=1}^n \psi_k^{(0)}(z)}}\leq
\Coeff{z^j}{\prod_{k=1}^n \psi_k^{(1)}(z)}\leq
\Coeff{z^j}{\prod_{k=1}^n \psi_k^{(2)}(z)}.&
\end{eqnarray*}
\end{lemma}
\Proofof{Lemma~\ref{la1a}}
We first note that 
\[\prod_{j=1}^nF_j=\Big(\prod_{j=1}^n(K_j-K_0+K_0)\Big)
\ee^{n\overline{L}}=\sum_{k=0}^nV_kK_0^{n-k}\ee^{n\overline{L}}
=\sum_{k=0}^n V_k\,G^{n-k}\ee^{k\overline{L}}=W_n.\]
For $k\in\set{n}$, let 
$\lambda(n,k)=0$ or $\lambda(n,k)=n-k-k\floor{(n-k)/k}$ 
according to whether $G$ is infinitely divisible or not. 
For $\ell\in\setn{n}$, we obtain
\begin{eqnarray} 
\prod_{j=1}^nF_j-W_\ell
&=&\sum_{k=\ell+1}^nV_kG^{n-k}
\ee^{k\overline{L}}
=\sum_{k=\ell+1}^n\sum_{J\subseteq{\set{n}}:\;\card{J}=k}
\Big(\prod_{j\in J}M_{j,k}\Big)G^{\lambda(n,k)}\nonumber\\
&=&\sum_{k=\ell+1}^n\coeff{z^k}{\psi_k(z)}G^{\lambda(n,k)},
\label{eq4454}
\end{eqnarray}
where $\psi_k(z)=\prod_{j=1}^n(\dirac+M_{j,k}z)$ 
is regarded as a formal power series for $k\in\set{n}$.
It should be mentioned that it is essential here to extract 
the~$k$th coefficient of a formal power series which itself depends 
on~$k$. By Lemma~\ref{la176}, for $k\in\set{n}$, we get
\begin{equation}\label{eq1876} 
\norm{\coeff{z^k}{\psi_k(z)}}
\leq\Coeff{z^k}{\prod_{j=1}^n (1+\norm{M_{j,k}}z)}\leq
\coeff{z^k}{\ee^{\nu_{k,1}z}}=\frac{\nu_{k,1}^k}{k!}.
\end{equation}
On the other hand, using
\begin{eqnarray*} 
\psi_k(z)&=&\Expo{\sum_{j=1}^nM_{j,k}z}
\prod_{j=1}^n\big(\ee^{-M_{j,k}z}(\dirac+M_{j,k}z)\big)\\
&=&\Expo{\sum_{j=1}^nM_{j,k}z}\prod_{j=1}^n
\Big[\sum_{m=0}^\infty\frac{1-m}{m!}(-M_{j,k})^mz^m\Big],
\end{eqnarray*}
we derive
\begin{equation}\label{eq5244}
\norm{\coeff{z^k}{\psi_k(z)}}\leq\Coeff{z^k}{\ee^{\widetilde{\nu}_kz}
\prod_{j=1}^n g(\norm{M_{j,k}}z)},
\end{equation}
where, for $y\in\CC$,
\[g(y)=\sum_{m=0}^\infty \frac{\ab{1-m}}{m!}y^m=
2-(1-y)\ee^y=1+\frac{y^2}{2}+\frac{y^3}{3}+\frac{y^4}{8}+\dots.\]
From the definition of $c_1^{}$, we obtain that
\begin{equation}\label{eq4521}
\ab{g(y)}\leq g(\ab{y})\leq \ee^{c_1\ab{y}^2}.
\end{equation}
Here, we note that $h(x)\defby\ln(2-(1-x)\ee^x)/x^2$ for 
$x\in(0,\infty)$ attains its maximum~$c_1^{}=0.694025\dots$ at 
point $x_0=0.936219\dots$. This can easily be shown using 
the representation
\[h(x)=\frac{1}{x^2}\int_0^1\frac{\dd}{\dd t}\ln(2-(1-tx)\ee^{tx})\,
\dd t=\int_0^1\frac{t}{2\ee^{-t x}-1+tx}\,\dd t,\]
which, after differentiation of the integrand, leads to a useful
integral formula of the derivative
\[h'(x)=\frac{1}{x^3}\int_0^x
\frac{t^2(2\ee^{-t}-1)}{(2\ee^{-t}-1+t)^2}\,\dd t.\]
As a consequence, we learn that $h'(x)=0$ has exactly one positive
solution $x=x_0$, which can easily be calculated numerically. Let
\[\bessel(0;y)=\sum_{m=0}^\infty\frac{(y^2/4)^m}{(m!)^2}=
\frac{1}{2\pii}\int_{-\pii}^\pii\expo{y \cos(t)}\,\dd t,
\qquad(y\in\CC)\]
be the modified Bessel function of first kind and order $0$. 
Using~(\ref{eq5244}), Cauchy's integral formula, and~(\ref{eq4521}), 
we derive, for $k\in\set{n}$ and arbitrary
$R_k\in(0,\infty)$, 
\begin{eqnarray*}
\norm{\coeff{z^k}{\psi_k(z)}}
&\leq&\frac{1}{2\pii\,R_k^k}\int_{-\pii}^{\pii} \ee^{-\ii k t}
\Big(\prod_{j=1}^n g(\norm{M_{j,k}}R_k\ee^{\ii t})\Big)
\expo{\widetilde{\nu}_k^{}R_k\ee^{\ii t}}\,\dd t\\
&\leq&\frac{1}{2\pii\,R_k^k}\int_{-\pii}^{\pii}
\expo{\widetilde{\nu}_k^{}R_k\cos(t)}\,\dd t
\prod_{j=1}^ng(\norm{M_{j,k}}R_k) \\
&\leq&\frac{1}{R_k^k}\bessel(0;\widetilde{\nu}_k^{}R_k)
\expo{c_1^{}\,\nu_{k,2}^{}\,R_k^2}\\
&=&\frac{\varphi(\widetilde{\nu}_k^{}R_k)}{R_k^k}
\Expo{\Big(\frac{\widetilde{\nu}_k^2}{4}\,
+c_1^{}\,\nu_{k,2}^{}\Big)R_k^2},
\end{eqnarray*}
where $\varphi(x)=\bessel(0;x)\,\ee^{-x^2/4}\leq1$, $(x\in\RR)$.
Choosing 
\[R_k=\Big(\frac{k}{2(4^{-1}\widetilde{\nu}_k^2
+c_1^{}\,\nu_{k,2}^{})}\Big)^{1/2},\]
we get 
\begin{equation}\label{eq597}
\norm{\coeff{z^k}{\psi_k(z)}}\leq\Big(\frac{2\ee\,c_1^{}}{k}\Big(
\frac{\widetilde{\nu}_k^2}{4c_1^{}}+\nu_{k,2}^{}\Big)\Big)^{k/2}.
\end{equation}
Taking into account~(\ref{eq4454}), the fact that 
$\lambda(n,k)\in\setn{n}$, as well as~(\ref{eq1876}) 
and~(\ref{eq597}), we obtain
\begin{eqnarray*} 
\Norm{\prod_{j=1}^nF_j-W_\ell}
&\leq&\sum_{k=\ell+1}^n\norm{\coeff{z^k}{\psi_k(z)}}
\leq\sum_{k=\ell+1}^n
\Big[\Big(\frac{2\ee\,c_1^{}}{k}
\Big(\frac{\widetilde{\nu}_k^2}{4c_1^{}}+\nu_{k,2}^{}\Big)
\Big)^{k/2}\wedge\frac{\nu_{k,1}^k}{k!}\Big].
\end{eqnarray*}
The proof is completed.~\hfill\qed
\subsection{Remaining proofs}\label{remproofs}
\Proofof{Theorem~\ref{thm01}}
Part~(a) follows from Lemma~\ref{la1a}.
Indeed, for $\eta_{\ell}^{}<(2\ee\,c_1^{})^{-1}$, we have
\begin{eqnarray*} 
\Norm{\prod_{j=1}^nF_j-W_\ell}&\leq&\sum_{k=\ell+1}^n
\Big(\frac{2\ee\,c_1^{}}{k^{1+\alpha}}
\Big(\frac{\widetilde{\nu}_k^2}{4c_1^{}}+\nu_{k,2}^{}\Big)\Big)^{
(\ell+1)/2}k^{\alpha(\ell+1)/2}\Big(\frac{2\ee\,c_1^{}}{k}
\Big(\frac{\widetilde{\nu}_k^2}{4c_1^{}}+\nu_{k,2}^{}\Big)\Big)^{
(k-\ell-1)/2}\\
&\leq&(2\ee\,c_1^{}\,\eta_{\ell,\alpha}^{})^{(\ell+1)/2}
\sum_{k=\ell+1}^nk^\beta(2\ee\,c_1^{}\,\eta_\ell^{})^{(k-\ell-1)/2}\\
&\leq&
(\ell+1)^{\beta}\,\beta!\,
\frac{(2\ee\,c_1^{}\,\eta_{\ell,\alpha}^{})^{(\ell+1)/2}}{
(1-\sqrt{2\ee\,c_1^{}\,\eta_{\ell}^{}})^{\beta+1}}.
\end{eqnarray*}
Here we used that, for $x\in[0,1)$, 
\[\sum_{k=\ell+1}^n k^\beta x^{k-\ell-1}
\leq(\ell+1)^\beta
\sum_{k=0}^\infty \frac{(k+\beta)!}{k!} x^k
=(\ell+1)^\beta\frac{\dd}{\dd x^\beta} \frac{1}{1-x}
=(\ell+1)^\beta\frac{\beta!}{(1-x)^{\beta+1}}.\]
Part~(b) is shown by using Corollary~\ref{cor227} together
with Remark~\ref{rem62}. In fact, for $j,k\in\set{n}$, we obtain
\[\norm{M_{j,k}}^2
\leq\Big(\norm{\restr{(F_j-G)}{B_j}}+\ab{(F_j-G)(B_j)}
+\sqrt{\frac{2k}{n}}\Big(\int_{B^\cc_j}(f_j-1)^2\,\dd G
\Big)^{1/2}\Big)^2,\]
since, for $k\in\set{n}$,
\[\Floor{\frac{n-k}{k}}+1
\geq \frac{\max\{n-k,k\}}{k}\geq\frac{n}{2k}.\]
Similarly,  we have
\[\widetilde{\nu}_k^2=\Norm{\sum_{j=1}^nM_{j,k}}^2
\leq \Big(n\norm{\restr{(\overline{F}-G)}{B_0}}+
n\ab{(\overline{F}-G)(B_0)}+\sqrt{2kn}\Big(\int_{B^\cc_0}
(f_0-1)^2\,\dd G\Big)^{1/2}\Big)^2.\]
This yields~(\ref{eq512}) and completes the 
proof.~\hfill\qed\medskip\\
For the proof of Theorem~\ref{thm02}, we need the following lemma. 
\begin{lemma} \label{la66}
Let $n\in\NN$, $L_1,\dots,L_n\in\M$ with $\sum_{j=1}^nL_j=0$,
and, for $k\in\setn{n}$ and $m\in\NN$,
\[\widetilde{V}_k=\sum_{J\subseteq\,\set{n}:\;\card{J}=k}
\prod_{j\in J}L_j,\qquad \vartheta_m=\sum_{j=1}^n\norm{L_j}^m.\]
Then we have 
\begin{eqnarray*}
&\displaystyle \norm{\widetilde{V}_2}\leq\frac{1}{2}\vartheta_{2},
\qquad
\norm{\widetilde{V}_3}\leq\frac{1}{3}\,\vartheta_{3},\qquad
\norm{\widetilde{V}_4}\leq\frac{1}{8}\,\vartheta_{2}^2,&\\
&\displaystyle\norm{\widetilde{V}_5}
\leq\frac{1}{6}\vartheta_{2}\vartheta_{3},\qquad
\norm{\widetilde{V}_6}\leq\frac{5}{144}\,\vartheta_{2}^3,\qquad
\norm{\widetilde{V}_7}\leq\frac{1}{24}\vartheta_{2}^2\vartheta_{3}.&
\end{eqnarray*}
\end{lemma}
\Proof 
The first two inequalities are easy.  Taking into account
(\ref{eq612})--(\ref{eq77}), it is not difficult to  show that, 
letting
$\widetilde{\Gamma}_m=\sum_{j=1}^n(-L_j)^m$,  $(m\in\NN)$,
\begin{eqnarray*}
\norm{\widetilde{V}_4}&=&\frac{1}{8}\norm{[\widetilde{\Gamma}_2^2
-\widetilde{\Gamma}_4]-\widetilde{\Gamma}_4}
\leq\frac{1}{8}\big(\norm{\widetilde{\Gamma}_2^2-\widetilde{\Gamma}_4}
+\norm{\widetilde{\Gamma}_4}\big)
\leq\frac{1}{8}\,\vartheta_{2}^2,\\
\norm{\widetilde{V}_5}&=&\frac{1}{6}\Norm{
[\widetilde{\Gamma}_2\widetilde{\Gamma}_3-\widetilde{\Gamma}_5]
-\frac{1}{5}\widetilde{\Gamma}_5}
\leq\frac{1}{6}\vartheta_{2}\vartheta_{3},\\
\norm{\widetilde{V}_6}&=&\frac{1}{144}
\norm{
3[\widetilde{\Gamma}_2^3-3\widetilde{\Gamma}_2\widetilde{\Gamma}_4
+2\widetilde{\Gamma}_6]
-9[\widetilde{\Gamma}_2\widetilde{\Gamma}_4-\widetilde{\Gamma}_6]
-8[\widetilde{\Gamma}_3^2-\widetilde{\Gamma}_6]
+\widetilde{\Gamma}_6}\\
&\leq&\frac{1}{144}\big(
3[\vartheta_{2}^3-3\vartheta_{2}\vartheta_{4}+2\vartheta_{6}]
+9[\vartheta_{2}\vartheta_{4}
-\vartheta_{6}]+8[\vartheta_{3}^{2}-\vartheta_{6}]
+\vartheta_{6}\big)\\
&=&\frac{1}{144}\big(
3\vartheta_{2}^3
+8[\vartheta_{3}^{2}-\vartheta_{6}]
-2\vartheta_{6}
\big)\leq\frac{5}{144}\,\vartheta_{2}^3,\\
\norm{\widetilde{V}_7}&=&
\frac{1}{840}\norm{
35[\widetilde{\Gamma}_2^2\widetilde{\Gamma}_3
-2\widetilde{\Gamma}_2\widetilde{\Gamma}_5
-\widetilde{\Gamma}_3\widetilde{\Gamma}_4+2\widetilde{\Gamma}_7]
-14[\widetilde{\Gamma}_2\widetilde{\Gamma}_5-\widetilde{\Gamma}_7]
-35[\widetilde{\Gamma}_3\widetilde{\Gamma}_4-\widetilde{\Gamma}_7]
+\widetilde{\Gamma}_7}\\
&\leq&\frac{1}{840}\big(
35[\vartheta_{2}^2\vartheta_{3}-2\vartheta_{2}\vartheta_{5}
-\vartheta_{3}\vartheta_{4}+2\vartheta_{7}]
+14[\vartheta_{2}\vartheta_{5}-\vartheta_{7}]
+35[\vartheta_{3}\vartheta_{4}-\vartheta_{7}]
+\vartheta_{7}
\big)\\
&=&\frac{1}{840}\big(28[\vartheta_{2}^2\vartheta_{3}
-2\vartheta_{2}\vartheta_{5}+\vartheta_{7}]
+7\vartheta_{2}^2\vartheta_{3}-6\vartheta_{7}\big)
\leq\frac{1}{24}\vartheta_{2}^2\vartheta_{3}.
\end{eqnarray*}
Observe that, in order to obtain good constants, a convenient grouping
of terms is essential. Further, for the bound of 
$\norm{\widetilde{V}_6}$, we used the inequality 
$(\vartheta_{3}^2-\vartheta_{6})/\vartheta_{2}^3\leq4^{-1}$, 
which can be proved by using
\[\frac{\vartheta_{3}^2-\vartheta_{6}}{\vartheta_{2}^3}
=\Big(\sum_{j=1}^nx_j^{3/2}\Big)^2-\sum_{j=1}^nx_j^3
\leq\sum_{j=1}^nx_j^2(1-x_j)\bydef g_n((x_1,\dots,x_n)),\]
where $x_j=\norm{L_j}^2(\sum_{i=1}^n\norm{L_i}^2)^{-1}$,
and the fact that the functions $g_n((x_1,\dots,x_n))$
for $n\in\{3,4,\dots\}$ and $(x_1,\dots,x_n)\in[0,\,1]^n$
with $\sum_{j=1}^nx_j=1$ satisfy
\[g_{n}((x_1,\dots,x_n))\leq g_{n-1}((x_1+x_2,x_3,\dots,x_n)),\]
whenever $0\leq x_1\leq \dots\leq x_n\leq 1$.
This completes the proof of the lemma.\hfill\qed\bigskip\\
\Proofof{Theorem~\ref{thm02}}
In order to prove the assertions, we need a further bound. 
In fact, similarly as in the proof of Lemma~\ref{la1a}, we get that,
for $\ell\in\setn{n}$,
\begin{eqnarray*}
\Norm{\prod_{j=1}^nF_j-W_\ell}
&\leq&1+\norm{W_\ell}\leq 
2+\Norm{\sum_{k=1}^{\ell}\Coeff{z^k}{\prod_{j=1}^n
\big(\dirac+M_{j,k}z\big)}G^{\lambda(n,k)}}\\
&\leq&2+\sum_{k=1}^\ell t^k=\frac{2-t-t^{\ell+1}}{1-t},
\end{eqnarray*}
where $t=\sqrt{2\ee\,c_1^{}\,\eta_0^{}}$. 
Similarly, if $G=\overline{F}$, then, for $\ell\in\set{n}$,
\begin{eqnarray*}
\Norm{\prod_{j=1}^nF_j-W_\ell}\leq2+\sum_{k=2}^\ell 
\widetilde{t}^k=\frac{2-2\widetilde{t}+\widetilde{t}^2
-\widetilde{t}^{\ell+1}}{1-\widetilde{t}},
\end{eqnarray*}
where $\widetilde{t}=\sqrt{2\ee\,c_1^{}\,\eta_1^{}}$. 
We now prove (a). Let $\ell\in\setn{n}$. If $t\in[0,x_\ell]$, 
then~(\ref{eq384b}) gives
\[\Norm{\prod_{j=1}^nF_j-W_\ell}\leq
\frac{t^{\ell+1}}{1-t}\leq\frac{t^{\ell+1}}{1-x_\ell}.\]
On the other hand, if $t\in(x_\ell,\,\infty)$, then 
\[\Norm{\prod_{j=1}^nF_j-W_\ell}\leq
\frac{2-t-t^{\ell+1}}{1-t}
=\frac{t^{\ell+1}}{1-x_\ell}
\Big(\frac{2-t-t^{\ell+1}}{t^{\ell+1}(1-t)}\Big)
\Big(\frac{2-x_\ell-x_\ell^{\ell+1}}{x_\ell^{\ell+1}(1-x_\ell)}
\Big)^{-1}\leq\frac{t^{\ell+1}}{1-x_\ell},\]
since $\frac{2-t-t^{\ell+1}}{t^{\ell+1}(1-t)}=t^{-\ell-1}
+\sum_{j=1}^{\ell+1}t^{-j}$ is decreasing on $(0,\infty)$.  This
yields~(\ref{eq385}) and~(\ref{eq3981}).  The proof of (a) is  easily
completed.  Let us now show (b). Set $G=\overline{F}$.  Similarly to
the above, one can show that, for $\ell\in\set{n}$, 
\[\Norm{\prod_{j=1}^nF_j-W_\ell}\leq\frac{\widetilde{t}^{\ell+1}}{
1-\widetilde{x}_\ell}.\]
This proves one part of (\ref{eq199}).
Using the norm inequalities in Lemma~\ref{la66} and~(\ref{eq384b}), 
we derive, for $\ell\in\set{3}$, 
\[\Norm{\prod_{j=1}^nF_j-W_\ell}\leq\Norm{\prod_{j=1}^nF_j-W_7}
+\Norm{W_7-W_\ell}\leq\Norm{\prod_{j=1}^nF_j-W_7}
+\sum_{k=\ell+1}^7 \norm{V_k\,\overline{F}^{n-k}}\leq 
\zeta_\ell(\eta_1^{}),\]
where, for $x\in[0,(2\ee\,c_1^{})^{-1})$, 
\begin{eqnarray*} 
\zeta_1(x)&=&x+\sqrt{3}\,x^{3/2}+2\,x^2+\frac{5^{5/2}}{6}\,x^{5/2}
+\frac{15}{2}\,x^3+\frac{7^{7/2}}{24}x^{7/2}+
\frac{(2\ee\,c_1^{}\,x)^{4}}{1-\sqrt{2\ee\,c_1^{}\,x}},\\
\zeta_2(x)&=&\sqrt{3}\,x^{3/2}+2\,x^2+\frac{5^{5/2}}{6}\,x^{5/2}
+\frac{15}{2}\,x^3+\frac{7^{7/2}}{24}x^{7/2}+
\frac{(2\ee\,c_1^{}\,x)^{4}}{1-\sqrt{2\ee\,c_1^{}\,x}},\\
\zeta_3(x)&=&2\,x^2+\frac{5^{5/2}}{6}\,x^{5/2}
+\frac{15}{2}\,x^3+\frac{7^{7/2}}{24}x^{7/2}+
\frac{(2\ee\,c_1^{}\,x)^{4}}{1-\sqrt{2\ee\,c_1^{}\,x}}.
\end{eqnarray*}
Note that, for $\ell\in\set{3}$, we have 
$\zeta_\ell(\eta_1^{})\leq\frac{2-2\widetilde{t}+\widetilde{t}^2
-\widetilde{t}^{\ell+1}}{1-\widetilde{t}}$, if and only if  
$\eta_1^{}\in[0,s_\ell]$, where $s_1=0.182839\dots$, 
$s_2=0.196439\dots$,
and $s_3=0.205094\dots$. If $\eta_1^{}\in[0,s_\ell]$, then
$\norm{\prod_{j=1}^nF_j-W_\ell}\leq \zeta_\ell(\eta_{1}^{})
\leq \zeta_\ell(s_\ell)\eta_{1}^{(\ell+1)/2}/s_\ell^{(\ell+1)/2}$. 
If $\eta_{1}^{}\in(s_\ell,\infty)$, then, letting 
$\widetilde{t}_\ell=\sqrt{2\ee\,c_1^{}s_\ell}$, 
\[\Norm{\prod_{j=1}^nF_j-W_\ell}
\leq\widetilde{t}^{\ell+1}\frac{2-2\widetilde{t}+\widetilde{t}^2
-\widetilde{t}^{\ell+1}}{\widetilde{t}^{\ell+1}(1-\widetilde{t})}
\leq\widetilde{t}^{\ell+1}\frac{2-2\widetilde{t}_\ell+
\widetilde{t}_\ell^2-\widetilde{t}_\ell^{\ell+1}}{
\widetilde{t}_\ell^{\ell+1}(1-\widetilde{t}_\ell)}
=\eta_1^{(\ell+1)/2}\frac{\zeta_\ell(s_\ell)}{s_\ell^{(\ell+1)/2}}.\]
Numerical calculations give the bounds for 
$\widetilde{u}_\ell$, $(\ell\in\set{3})$ as claimed in (\ref{eq199}).
This completes the proof.~\hfill\qed\bigskip\\
\Proofof{Proposition~\ref{prop256}} Consider fixed $j,k\in\set{n}$.
Let $d=2b$, $\overline{p}=(\overline{p}_1,\dots,\overline{p}_b,
\overline{p}_1,\dots,\overline{p}_b)\in\Zpl^d$, and
$\rho=\floor{(n-k)/k}$. Further, for
$v\in\Zpl^{d}$ with $\vecsum{v}=2$, 
let $a_v=\overline{p}_r-p_{j,r}$ if 
$v_r=v_{b+r}=1$ and $a_v=0$ otherwise.
Let
\[H_0=\dirac,\qquad H_r=\left\{\begin{array}{ll}
\dirac_{-x_r},& r\in\set{b},\\
\dirac_{x_{r-b}},& r\in\set{d}\setminus\set{b}.
\end{array}\right.\]
Then we have $\overline{F}=\sum_{r=0}^d\overline{p}_rH_r$ and
\begin{eqnarray*} 
F_j-\overline{F}&=&(p_{j,0}-\overline{p}_0)\dirac+
\sum_{r=1}^b(p_{j,r}-\overline{p}_r)
(\dirac_{-x_r}+\dirac_{x_r})\\
&=&\sum_{r=1}^b(p_{j,r}-\overline{p}_r)
(\dirac_{-x_r}+\dirac_{x_r}-2\dirac)
=\sum_{r=1}^b(\overline{p}_r-p_{j,r})
(\dirac_{-x_r}-\dirac)(\dirac_{x_r}-\dirac)\\
&=&\sum_{\vecsum{v}=2}\frac{a_v}{v!}\prod_{r=1}^d(H_r-H_0)^{v_r}.
\end{eqnarray*}
Here and henceforth, sums over $v$ and $w$ are taken over subsets
of  $\Zpl^d$ as indicated. In particular, we obtain
\[\norm{M_{j,k}}^2
=\norm{(F_j-\overline{F})\overline{F}^{\rho}}^2
\leq\norm{F_j-\overline{F}}^2
\leq \Big(\ab{p_{j,0}-\overline{p}_0}+
2\sum_{r=1}^b\ab{p_{j,r}-\overline{p}_r}\Big)^2.\]
On the other hand, in view of
\[\norm{M_{j,k}}^2=\norm{(F_j-\overline{F})\overline{F}^{\rho}}^2
=\Norm{\Big(\sum_{\vecsum{v}=2}\frac{a_v}{v!}
\prod_{r=1}^d(H_r-H_0)^{v_r}\Big)
\Big(\sum_{r=0}^d\overline{p}_rH_r\Big)^\rho}^2,\]
we see that (\ref{eq551a}) can be applied, which together with the  
simple fact that 
$\frac{\rho!}{(\rho+2)!}\leq\frac{1}{(\rho+1)^2}\leq\frac{4k^2}{n^2}$
gives
\begin{eqnarray*} 
\norm{M_{j,k}}^2&\leq&\frac{\rho!}{(\rho+2)!}
\sum_{\vecsum{w}\leq 2}\frac{w!(2-\vecsum{w})!}{
\overline{p}^w\,\overline{p}_0^{2-\vecsum{w}}}
\Big[\sum_{\vecsum{v}=2}\frac{a_{v}}{v!}
\prod_{r=1}^{d}
\binomial{v_r}{w_r}\Big]^2\\
&\leq&\frac{4k^2}{n^2}\Big(\frac{2}{\overline{p}_0^{2}}
\Big[\sum_{\vecsum{v}=2}\frac{a_{v}}{v!}\Big]^2 
+\sum_{\vecsum{w}=1}\frac{1}{\overline{p}^w\,
\overline{p}_0}
\Big[\sum_{\vecsum{v}=2}\frac{a_{v}}{v!}
\prod_{r=1}^{d}
\binomial{v_r}{w_r}\Big]^2\\
&&\hspace{5cm}{}
+\sum_{\vecsum{w}=2}\frac{w!}{\overline{p}^w}
\Big[\sum_{\vecsum{v}=2}\frac{a_{v}}{v!}
\prod_{r=1}^{d}
\binomial{v_r}{w_r}\Big]^2\Big).
\end{eqnarray*}
The special definition of $a_v$, $(v\in\Zpl^d,\vecsum{v}=2)$
implies that $a_{\unitvec{r(1)}+\unitvec{r(2)}}=0$ for $r(1),r(2)\in\set{d}$
with $\ab{r(1)-r(2)}\neq b$ and therefore
the terms on the right-hand side can be evaluated as follows:
\begin{eqnarray*}
\frac{2}{\overline{p}_0^2}
\Big[\sum_{\vecsum{v}=2}\frac{a_{v}}{v!}\Big]^2&=&
\frac{2}{\overline{p}_0^2}
\Big[\sum_{r=1}^b(\overline{p}_r-p_{j,r})\Big]^2
=\frac{(p_{j,0}-\overline{p}_0)^2}{2\overline{p}_0^2},\\
\sum_{\vecsum{w}=1}\frac{1}{\overline{p}^w\,
\overline{p}_0}
\Big[\sum_{\vecsum{v}=2}\frac{a_{v}}{v!}
\prod_{r=1}^{d}\binomial{v_r}{w_r}\Big]^2&=&
\sum_{r=1}^d\frac{1}{\overline{p}_r\,\overline{p}_0}
\Big[\sum_{\vecsum{v}=2}\frac{a_{v}}{v!}
\binomial{v_r}{1}\Big]^2
=2\sum_{r=1}^b\frac{(\overline{p}_r-p_{j,r})^2}{
\overline{p}_r\,\overline{p}_0},\\
\sum_{\vecsum{w}=2}\frac{w!}{\overline{p}^w}
\Big[\sum_{\vecsum{v}=2}\frac{a_{v}}{v!}
\prod_{r=1}^{d}
\binomial{v_r}{w_r}\Big]^2&=&
\sum_{\vecsum{w}=2}\frac{w!}{\overline{p}^w}
\Big[\frac{a_{w}}{w!}\Big]^2=
\sum_{r=1}^b\frac{(\overline{p}_r-p_{j,r})^2}{
\overline{p}_r^2}.
\end{eqnarray*}
We note that some of the binomial coefficients above are equal to
zero. This implies that
\begin{eqnarray*}
\norm{M_{j,k}}^2
&\leq&\frac{4k^2}{n^2}\Big(\frac{(\overline{p}_0-p_{j,0})^2}{
2\overline{p}_0^{2}}+2\sum_{r=1}^b
\frac{(\overline{p}_r-p_{j,r})^2}{\overline{p}_r
\overline{p}_0}+\sum_{r=1}^b
\frac{(\overline{p}_r-p_{j,r})^2}{\overline{p}_r^2}\Big).
\end{eqnarray*}
Using this together with  
\[\eta_{\ell,1}^{}=\max_{k\,\in\,\set{n}\setminus\set{\ell}}
\frac{\nu_{k,2}^{}}{k^{2}}=\max_{k\,\in\,\set{n}\setminus\set{\ell}}
\Big(\frac{1}{k^{2}}\sum_{j=1}^n\norm{M_{j,k}}^2\Big),\]
(see the comment after Theorem~\ref{thm01})
the proof is easily completed.~\hfill\qed
\section*{Acknowledgment}
The author thanks Lutz Mattner for helpful discussions. 
He is also grateful to both referees for valuable remarks,
which led to an improvement of the paper.
\linespread{1.3}


\end{document}